\documentclass[a4paper, 10pt, conference]{ieeeconf}      

\IEEEoverridecommandlockouts                              

\usepackage{hyperref}
\usepackage{amsmath, amssymb, mathrsfs}
\usepackage{graphicx,graphics,color}
\usepackage{lineno}
\usepackage{ulem}
\usepackage{amssymb}

\newtheorem{theorem}{Theorem}[section]

\newtheorem{definition}{Definition}[section]

\title{\LARGE \bf
Traffic Modeling and Real-time Control for Metro Lines
}

\author{Nadir Farhi$^{*}$, Cyril Nguyen Van Phu, Habib Haj-Salem \& Jean-Patrick Lebacque 
\thanks{Universit\'e Paris-Est, Ifsttar / Cosys / Grettia, F-77447 Paris - Marne-la-Vall\'ee cedex 2, France.}
\thanks{$^{*}$ Corresponding author: {\tt\small nadir.farhi@ifsttar.fr}} } %

\begin{document}

\maketitle
\thispagestyle{empty}
\pagestyle{empty}

\begin{abstract}
We present in this article traffic flow and control models for the train dynamics in metro lines.
The first model, written in the max-plus algebra, takes into account minimum running, dwell and
safety time constraints, without any control of the train dwell times at platforms, and
without consideration of the passenger travel demand.
We show that the dynamics are stable and converge to stationary regimes with a unique asymptotic average growth rate.
Moreover, the asymptotic average train time-headway, dwell time, as well as close-in time, are derived analytically, as
functions of the number of running trains on the metro line.
We then introduce, in a second model, the effect of the passenger demand on the train dwell times at platforms.
We review that, if this effect is not well controlled, then the traffic is unstable.
Finally, we propose a traffic control model which deals with this instability, by well controlling the effect of
passenger arrivals on the train dwell times at platforms.
We show that the dynamics are stable and converge to stationary regimes with a unique asymptotic average growth rate.
We then calculate by numerical simulations the asymptotic average time-headway as a function of the number of running trains,
compare the results with those of the max-plus algebra model, and derive the effect of the passenger travel demand on the
frequency of the metro line, under the proposed control model.
\end{abstract}

\textbf{Keywords.} Traffic control, traffic modeling, max-plus algebra, railway traffic.

\section{INTRODUCTION}
\label{introduction}

Although the mass transit by metro is known to be one of the most efficient, safe, and comfortable passenger transport mode
in urban and sub-urban areas, it is also known to be naturally unstable when exploited at high frequencies~\cite{BCB91}.
Indeed, at a high frequency, where the capacity margins are reduced, train delays are amplified and rapidly propagated in space
and in time.
In severe perturbations on a given metro line, the delays may propagate to other lines.
One of the parameters affecting the dwell times of trains at platforms, and by that causing train delays, is the passenger flow time.
It is the largest and the hardest parameter to control~\cite{TCQSM}.
If no efficient traffic control is set for the metro line, a small delay may rapidly drive down the train frequency.
In order to deal with that, innovative approaches and methods for real-time railway traffic control are needed.
We propose in this article discrete-event traffic models, with a
control strategy that guarantees train dynamics stability, and that
takes into account the effect of passenger arrivals. 

We are concerned here by real-time control of the train dynamics in metro lines.
As well known, one of the important control parameters of the train dynamics is the train dwell
times at platforms. Indeed, on inter-station links, trains follow generally fixed speed and acceleration profiles
which they adapt to the train dynamics with respect to their respective leading trains 
(directly in case of a moving block system, or indirectly by means of traffic lights, 
in case of a fixed block system).
The control of the train dynamics is not trivial because of the passenger effect on the train dwell times 
at platforms. In effect, high passenger densities on platforms and/or in trains induce additional constraints for the train dwell times at platforms.
This situation can be caused by high level of passenger demand, or by one or more than one delayed train(s) (incidents).
A direct consequence of those constraints is the extension of the train dwell times at platforms.
The caused train delays then propagate in space and in time and extend the passenger waiting and travel times.

Real-time railway traffic control has been treated since decades, with several approaches (mathematical, simulation-based, expert system, etc.)
Cury et al. (1980)~\cite{Cur80} proposed an analytic traffic model with a multilevel hierarchical optimization method.
Breusegem et al. (1991)~\cite{BCB91} developed discrete event traffic and control models, pointing out the nature of traffic instability on metro lines.
The authors of~\cite{BCB91} proposed a linear quadratic (LQ) control approach to deal with the instability.
Lee et al. (1997)~\cite{Lee97} derived by simulation minimum time-headways on the orange line of the Kaohsiung Mass Rapid Transit (KMRT),
and compared them to existing theoretical formulas.
Assis and Milani (2004)~\cite{Ass04} solved the train scheduling problem in metro lines by a predictive control method.
A stochastic model for the propagation of train delays across the network has been proposed by Engelhardt-Funke and Kolonko (2004)~\cite{EK04},
where the waiting time for passengers changing lines is minimized.

Goverd (2007)~\cite{Gov07} developed an algebraic model for the analysis of timetable stability and robustness against delays.
The approach permits to evaluate the realizability and the stability of timetables by max-plus spectral analysis,
and to assess their robustness against delays, using critical path algorithms.
The model has been applied to the Dutch national railway timetable.
A passenger boarding time model, where the train dwell time is a function of the train time-headway and of the passenger travel demand 
origin-destination (OD) matrix, is developed by Cao et al. (2009)~\cite{Cao09}.
Andre~\cite{And} studied the linear effect of the passenger alighting and boarding flows, and of the on-board crowding level,
on the train dwell time.
Sun et al. (2013)~\cite{SZDZ13} have proposed a mathematical programming approach to minimize waiting time, 
where the Lagrangian duality theory is used to solve the optimization problem.
Train schedules that take into account passengers' exchanges at platforms are generated.
More recently, Cacchiani et al. (2014)~\cite{Cac14} have given an overview of recovery models and algorithms for real-time 
railway disturbance and disruption management.

The problem for which we propose a solution here is the following.
In one side, the train dwell times at platforms need to be extended to respond to an increasing of passenger densities at platforms,
or of passenger arrivals.
In the other side, by adopting such a law, the accumulation of passengers at a platform, caused by a delayed train arrival,
would extend the dwell time of the train at the considered platform, 
which would induce more delay for the same train for its arrival to the next platform,
and which would delay the departures from the upstream stations of the considered station.
Therefore, delays are propagated and amplified, and the traffic is unstable.
The traffic control model we present here proposes a good compromise for the control of train dwell times at platforms,
that satisfies passengers travel demand, and that assures stable train dynamics.

The outline of this article is as follows.
In Section~2, we give a short review on the dynamic programming systems.
In Section~3, we present a max-plus algebra discrete-event dynamic model that does not take into account the passenger arrivals to platforms.
The model permits, in particular, an analytical derivation of traffic phase diagrams for the train dynamics.
In Section~4, we extend the max-plus model in order to take into account the effect of the passenger arrivals on the train dynamics.
We briefly review the natural instability of the train dynamics in the case where the effect of the passenger demand is uncontrolled.
In Section~5, we propose a modification of the latter model in order to guarantee the stability of the train dynamics, in addition
to take into account the passenger arrivals. We show that the traffic dynamics are interpreted
as dynamic programming systems of stochastic optimal control problems of Markov chains. 
We give the whole interpretation of the parameters, and derive the characteristics of the dynamics.
The effect of the passenger arrival rates on the asymptotic average train time-headway is shown under
the proposed control strategy. Finally, we draw some conclusions.

\section{A short review on dynamic programming systems}
\label{sec-review}

We give in this section a short review on dynamic programming systems and associated optimal
control problems. We present the review in three subsections: A - General dynamic programming systems,
B - Dynamic programming systems of stochastic optimal control problems, and C - Dynamic programming systems
of deterministic optimal control problems (Max-plus linear algebra systems).

\subsection{General dynamic programming systems}
\label{sub-sec-gdps}

A map $\mathbf{f} : \mathbb R^n \to \mathbb R^n$ is said to be additive 1-homogeneous if it satisfies:
$\forall x \in \mathbb R^n, \forall a \in \mathbb R, \mathbf{f}(a \textbf{1} + x) = a\textbf{1} + \mathbf{f}(x)$,
where $\textbf{1} \stackrel{def}{=} {}^t(1, 1, \ldots , 1)$. 
It is said to be monotone if it satisfies
$\forall x, y \in \mathbb R^n, x \leq y \Rightarrow \mathbf{f}(x) \leq \mathbf{f}(y)$, where $x \leq y$ means $x_i \leq y_i \forall i, 1 \leq i \leq n$.
If $\mathbf{f}$ is 1-homogeneous and monotone, then it is non expansive (or 1-Lipschitz) for the supremum norm~\cite{CT80}, i. e.
$\forall x, y \in \mathbb R^n, ||\mathbf{f}(x) - \mathbf{f}(y)||_{\infty} \leq ||x - y||_{\infty}$.
In this case, a directed graph $\mathcal G(\mathbf{f})$ is associated to $\mathbf{f}$.

The directed graph $\mathcal G(\mathbf{f})$ associated to $\mathbf{f}: \mathbb R^n \to \mathbb R^n$~\cite{GG99} is defined by the
set of nodes $\{1, 2, \ldots,n\}$ and by a set of arcs such that there exists an arc from
a node $i$ to a node $j$ if $\lim_{\eta\to\infty} \mathbf{f}_i(\eta e_j ) = \infty$, where $e_j$ is the $j$th
vector of the canonical basis of $\mathbb R^n$. We review below an important result on the existence
of additive eigenvalues of 1-homogeneous monotone maps.

\begin{theorem}\cite{GG98b, GK95}\label{th-dps}
  If $\mathbf{f} : \mathbb R^n \to \mathbb R^n$ is 1-homogeneous and monotone and if $\mathcal G(\mathbf{f})$
  is strongly connected then $\mathbf{f}$ admits an (additive) eigenvalue, i.e.
  $\exists \mu \in \mathbb R, \exists x \in \mathbb R^n : \mathbf{f}(x) = \mu + x$.
  Moreover, we have $\chi (\mathbf{f}) = \mu \textbf{1}$, where $\chi (\mathbf{f})$ denotes the asymptotic average growth vector of the
  dynamic system $x(k+1) = \mathbf{f}(x(k))$, defined by: $\chi (\mathbf{f}) = \lim_{k \to \infty} \mathbf{f}^k(x)/k$.
\end{theorem}

For simplicity, we use in this article, for all the dynamic systems, the notation $x^k$ instead of $x(k)$.
We give in the following, a natural extension of Theorem~\ref{th-dps}, which will permit us to
consider dynamic systems of the form $x^k = \mathbf{f}(x^k,x^{k-1},\ldots,x^{k-m+1})$.
For that, let us consider $\mathbf{f}: \mathbb R^{m\times n} \to \mathbb R^n$
associating for  $(x^{(0)},x^{(1)},\cdots, x^{(m-1)})$, where $x^{(i)}$ are vectors in $\mathbb R^n$,
a column vector in $\mathbb R^n$.
$$\begin{array}{llll}
  \mathbf{f}: & \mathbb R^{m\times n} & \to & \mathbb R^n \\
     & (x^{(0)},x^{(1)},\cdots, x^{(m-1)}) & \mapsto & \mathbf{f}(x).
\end{array}$$
We denote by $\mathbf{f}_i, i=0,1,\ldots,m-1$, the following maps.
$$\begin{array}{llll}
  \mathbf{f}_i: & \mathbb R^{n} & \to & \mathbb R^n \\
     & x & \mapsto & \mathbf{f}_i(x) = \mathbf{f}(-\infty, \ldots,-\infty,x,-\infty, \ldots, -\infty).\\
     &   &         & \qquad \qquad \qquad \qquad \qquad \; \; \; \uparrow \\
     &   &         & \quad \qquad \qquad \qquad \qquad \; \; \; i^{\text{th}} \text{component} 
\end{array}$$
and by $\tilde{\mathbf{f}}$ the following map.
$$\begin{array}{llll}
  \tilde{\mathbf{f}}: & \mathbb R^{n} & \to & \mathbb R^n \\
             & x & \mapsto & \tilde{\mathbf{f}}(x) = \mathbf{f}(x,x, \ldots,x).
\end{array}$$

\begin{theorem}\label{th-dps2}
   If $\tilde{\mathbf{f}}$ is additive 1-homogeneous and monotone, and if $\mathcal G(\tilde{\mathbf{f}})$
   id strongly connected and $\mathcal G(\mathbf{f}_0)$ is acyclic, then $\mathbf{f}$ admits a unique generalized
   additive eigenvalue $\mu > -\infty$ and an additive eigenvector $v > -\infty$, such that
   $\mathbf{f}(v, v-\mu, v-2\mu, \ldots, v-(m-1)\mu) = v$.
   Moreover, $\chi(\mathbf{f}) = \mu \textbf{1}$.
\end{theorem}

\proof
The proof consists in showing that the dynamic system $x^k = \mathbf{f}(x^k,x^{k-1},\ldots,x^{k-m+1})$
is equivalent to another dynamic system $z^k = \mathbf{h}(z^{k-1})$, where $\mathbf{h}$ is built from $\mathbf{f}$,
such that $\mathbf{h}$ satisfies additive 1-homogeneity, monotonicity and connectivity properties needed by Theorem~\ref{th-dps}.
We give here a sketch of the proof. The whole proof is available in Appendix~\ref{prf-dps2}.
Two steps are needed for the proof.
\begin{enumerate}
  \item Eliminate the dependence of $x^k$ on $x^{k-2}$, $x^{k-3}$, $\ldots$, $x^{k-m+1}$. 
    This is done by the well known state augmentation technique.
  \item Eliminate implicit terms (the dependence of $x^k$ on $x^k$).
    This is done by defining an order of updating the $n$ components op $x^k$, in
    such a way that no implicit term appears. This is possible because $\mathcal G(\mathbf{f}_0)$ 
    is acyclic.
\end{enumerate}
\endproof

\subsection{Dynamic programming systems of stochastic optimal control problems}
\label{sub-sec-dps}

We review here a particular additive 1-homogeneous and monotone dynamic programming system encountered
in stochastic optimal control of Markov chains. The dynamic system is written
\begin{equation}\label{eq-dp2}
	  x_i^{k+1} = \max_{u\in\mathcal U} ([M^{u} x^k]_i + c^{u}_i), \quad \forall 1 \leq i \leq n.
\end{equation}
where $\mathcal U$ is a set of indexes corresponding to controls; $M^{u}$,
for $u\in \mathcal U$, are stochastic matrices (i.e. satisfy $M^{u}_{ij} \geq 0$ and $M^{u} \textbf{1} = \textbf{1}$);
and $c^{u}$, for $u\in \mathcal U$, are reward vectors in $\mathbb R^n$.

(\ref{eq-dp2}) is the dynamic programming system associated to the stochastic optimal control of a Markov chain
with state space $\{1,2,\ldots,n\}$, transition matrices $M^u, u\in\mathcal U$, and associated rewards $c^u, u\in\mathcal U$.
The variable $x_i^k, i=1,2,...,n, k\in\mathbb N$ is interpreted, in this case, as the function value
of the stochastic control problem.	

In the traffic models of the train dynamics we propose in this article, we are concerned with system~(\ref{eq-dp2}) above.
More precisely, we will consider dynamic systems with implicit terms~(\ref{eq-dp3}).
\begin{equation}\label{eq-dp3}
  x_i^{k} = \max_{u\in\mathcal U} ([M^{u} x^{k-1}]_i + [N^{u} x^k]_i + c^{u}_i), \quad \forall 1 \leq i \leq n,
\end{equation}
where $M^u, u\in\mathcal U$ and $N^u, u\in\mathcal U$ satisfy $M^u_{ij} \geq 0, N^u_{ij} \geq 0$, 
$\sum_j (M^u_{ij}) \leq 1$, $\sum_j (N^u_{ij}) \leq 1$, and $\sum_j (M^u_{ij}+N^u_{ij}) = 1, \forall i,u$. 
For the analysis of dynamic system~(\ref{eq-dp3}) we will use Theorem~\ref{th-dps2}.

We notice here that, in our traffic models, the interpretation of systems~(\ref{eq-dp2}) and~(\ref{eq-dp3}) will be different 
from the stochastic optimal control one, in the sense that, in our models, such systems model directly
the train dynamics, and are not derived from stochastic optimal control problems.

\subsection{Max-plus algebra}

As mentioned above, dynamic programming systems associated to deterministic optimal control problems can be written linearly in the Max-plus algebra.
We present here some reviews in this algebra. 
The first traffic model we propose in section~\ref{sec-mp}, is written in the Max-plus algebra of square matrices of polynomials.
We review here the construction of this algebraic structure, and give some results that we used in the analysis
of our models.

\subsubsection{Max-plus algebra of scalars ($\mathbb R_{\max}$)}

Max-plus algebra~\cite{BCOQ92} is the idempotent commutative semi-ring $(\mathbb R \cup \{-\infty\}, \oplus, \otimes)$
denoted by $\mathbb R_{\max}$, where the operations $\oplus$ and $\otimes$ are defined by:
$a \oplus b = \max\{a, b\}$ and $a \otimes b = a + b$. The zero element is $(-\infty)$ denoted by $\varepsilon$ and the
unity element is $0$ denoted by $e$. 

\subsubsection{Max-plus algebra of square matrices ($\mathbb R_{\max}^{n\times n}$)}

We have the same structure on the set of square matrices.
If $A$ and $B$ are two Max-plus matrices of size $n \times n$, the addition
$\oplus$ and the product $\otimes$ are defined by: $(A \oplus B)_{ij} = A_{ij} \oplus B_{ij} , \forall i, j$, and
$(A \otimes B)_{ij} = \bigoplus_k[A_{ik} \otimes B_{kj}]$. The zero and the unity matrices are still denoted
by $\varepsilon$ and $e$ respectively.

A matrix $A$ is said to be \textit{reducible} if there exists a permutation matrix $P$ such that $P^T A P$ is lower block triangular.
A matrix that is not reducible is said to be \textit{irreducible}.

For a matrix $A\in \mathbb R_{\max}^{n\times n}$, a precedence graph $\mathcal G(A)$ is associated.
The set of nodes of $\mathcal G(A)$ is $\{1,2,\ldots,n\}$. There is an arc from node $i$ to node $j$
in $\mathcal G(A)$ if $A_{ji} \neq \varepsilon$.
A graph is said to be strongly connected if there exists a path from any node to any other node.
$A\in \mathbb R_{\max}^{n\times n}$ is irreducible if and only if $\mathcal G(A)$ is strongly connected~\cite{BCOQ92}.
That is, $A$ is irreducible if $\forall 1\leq i,j \leq n, \exists m\in\mathbb N, (A^m)_{ij} \neq \varepsilon$.

\subsubsection{Max-plus algebra of polynomials ($\mathbb R_{\max}[X]$)}

A (formal) polynomial in an indeterminate $X$ over $\mathbb R_{\max}$ is a finite sum $\bigoplus_{l=0}^p a_l X^l$ for some integer $p$
and coefficients $a_l\in \mathbb R_{\max}$. The set of formal polynomials in $\mathbb R_{\max}$ is denoted $\mathbb R_{\max}[X]$.
The support of a polynomial $\mathbf{f} = \bigoplus_{l=0}^p a_l X^l$ is $Supp(\mathbf{f}) = \{l, 0\leq l\leq p, a_l > \varepsilon\}$.
The degree of $\mathbf{f}$ is $Deg(\mathbf{f}) = \bigoplus_{l\in Supp(\mathbf{f})} l$.
The addition and the product of two polynomials $\mathbf{f}=\bigoplus_{l=0}^p a_l X^l$ and $\mathbf{g}=\bigoplus_{l=0}^q b_l X^l$ 
in $\mathbb R_{\max}[X]$ are defined as follows.
$$\mathbf{f}\oplus \mathbf{g} := \bigoplus_{l=0}^{deg(\mathbf{f})\oplus deg(\mathbf{g})} (a_l\oplus b_l) X^l,$$
$$\mathbf{f} \otimes \mathbf{g} := \bigoplus_{l=0}^{def(\mathbf{f})\otimes deg(\mathbf{g})} \left( \bigoplus_{i\otimes j=l} a_i\otimes b_j \right) X^l.$$
The zero element is $\varepsilon = \varepsilon X^0$ and the unity element is $e = eX^0$.
We notice that $\forall x\in \mathbb R_{\max}$, the valuation mapping $\varphi: \mathbf{f} \mapsto \mathbf{f}(x)$ is a homomorphism
from $\mathbb R_{\max}[X]$ into $\mathbb R_{\max}$.

\subsubsection{Max-plus algebra of polynomial square matrices $\left(\mathbb R_{\max}[X]\right)^{n\times n}$}

A polynomial matrix $A(X) \in \left(\mathbb R_{\max}[X]\right)^{n\times n}$ is a matrix with polynomial entries
$A_{ij}(X) = \bigoplus_{l=0}^p a_{ij}^{(l)} X^l$, where $a_{ij}^{(l)} \in \mathbb R_{\max}, \forall i,j, 0\leq i,j\leq n$ and
$\forall l, 1\leq l\leq p$.
$\left(\mathbb R_{\max}[X]\right)^{n\times n}$ is an idempotent semiring. The addition and the product are defined as follows.
$$(A(X) \oplus B(X))_{ij} = A_{ij}(X) \oplus B_{ij}(X), \forall i, j,$$
$$(A(X) \otimes B(X))_{ij} = \bigoplus_k[A_{ik}(X) \otimes B_{kj}(X)], \forall i, j.$$
The zero and the unity matrices are still denoted by $\varepsilon$ and $e$ respectively.
We also have, $\forall x\in \mathbb R_{\max}$, the valuation mapping $\varphi: A(X) \mapsto A(x)$ is a homomorphism 
from $\left(\mathbb R_{\max}[X]\right)^{n\times n}$ into $\mathbb R_{\max}^{n\times n}$.

A polynomial matrix $A (X) \in \left(\mathbb R_{\max}[X]\right)^{n\times n}$ is said to be \textit{reducible} if there exists a permutation
matrix $P(X)\in \left(\mathbb R_{\max}[X]\right)^{n\times n}$ such that $P(X)^T A(X) P(X)$ is lower block triangular.

Irreducibility of a matrix depends only on its support, that is the pattern of nonzero entries of the matrix.
As for finite values of $X$, the support of a matrix is preserved by the homomorphism valuation map $\varphi$,
then $A (X) \in \left(\mathbb R_{\max}[X]\right)^{n\times n}$ is irreducible if and only if $A(e)\in \left(\mathbb R_{\max}\right)^{n\times n}$ is so.
Therefore, $A(X) \in \left(\mathbb R_{\max}[X]\right)^{n\times n}$ is irreducible if and only if $\mathcal G(A(e))$ is strongly connected.

For a polynomial matrix $A(X) \in \left(\mathbb R_{\max}[X]\right)^{n\times n}$, a precedence graph $\mathcal G(A(X))$ is associated.
As for graphs associated to square matrices in $\mathbb R_{\max}^{n\times n}$, the set of nodes of $\mathcal G(A(X))$ is
$\{1,2,\ldots,n\}$. There is an arc $(i,j,l)$ from node $i$ to node $j$ in $\mathcal G(A(X))$ if $\exists l, 0\leq l\leq p, a^{(l)}_{ji} \neq \varepsilon$.
Moreover, a \textit{weight} $W(i,j,l)$ and a \textit{duration} $D(i,j,l)$ are associated to every arc $(i,j,l)$ in the graph,
with $W(i,j,l) = (A_l)_{ij} \neq \varepsilon$ and $D(i,j,l) = l$. Similarly, a weight, resp. duration of a cycle (directed cycle) in the graph is
the standard sum of the weights, resp. durations of all the arcs of the cycle. 
Finally, the \textit{cycle mean} of a cycle $c$ with a weight $W(c)$ and a duration $D(c)$ is $W(c)/D(c)$.

\subsubsection{Homogeneous linear Max-plus algebra systems}

We are interested in the first model we propose in this article, in the dynamics of a homogeneous $p$-order max-plus system
\begin{equation}\label{eq-mp1}
  x(k) = \bigoplus_{l=0}^p A_l \otimes x(k-l),
\end{equation}
where $x(k), k\in \mathbb Z$ is a sequence of vectors in $\mathbb R_{\max}^n$, and $A_l, 0\leq l\leq p$ are matrices in $\mathbb R_{\max}^{n\times n}$.
If we define $\gamma$ as the back-shift operator applied on the sequences of vectors  in $\mathbb Z$, such that: $\gamma x(k) = x(k-1)$,
and then more generally $\gamma^l x(k) = x(k-l), \forall l\in\mathbb N$, then~(\ref{eq-mp1}) is written as follows.
\begin{equation}\label{eq-mp2}
  x^k = \bigoplus_{l=0}^p \gamma^l A_l x^k = A(\gamma) x^k,
\end{equation}
where $A(\gamma)=\bigoplus_{l=0}^p \gamma^l A_l \in \left(\mathbb R_{\max}[\gamma]\right)^{n\times n}$ is a polynomial matrix in the back-shift
operator $\gamma$; see~\cite{BCOQ92,Gov07} for more details.

\begin{definition}
    $\mu \in \mathbb R_{\max} \setminus \{\varepsilon\}$ is said to be a \textit{generalized} eigenvalue~\cite{CCGMQ98} of $A(\gamma)$, with associated
    \textit{generalized} eigenvector $v\in \mathbb R_{\max}^n \setminus \{\varepsilon\}$, if $A(\mu^{-1}) \otimes v = v$,
    where $A(\mu^{-1})$ is the matrix obtained by evaluating the polynomial matrix $A(\gamma)$ at $\mu^{-1}$.
\end{definition}    

\begin{theorem} \cite[Theorem 3.28]{BCOQ92} \cite[Theorems 7.4.1 and 7.4.7]{Gov05} \label{th-mpa}
  Let $A(\gamma) = \oplus_{l=0}^p A_l\gamma^l$
  be an irreducible polynomial matrix with acyclic
  sub-graph $\mathcal G(A_0)$. Then $A(\gamma)$ has a unique generalized eigenvalue $\mu > \varepsilon$ and finite eigenvectors $v > \varepsilon$
  such that $A(\mu^{-1}) \otimes v = v$, and $\mu$ is equal to the maximum cycle mean of $\mathcal G(A(\gamma))$, given as follows.
  $\mu = \max_{c\in\mathcal C} W(c) / D(c)$,
  where $\mathcal C$ is the set of all elementary cycles in $\mathcal G(A(\gamma))$.
  Moreover, the dynamic system $x^k = A(\gamma) x^k$ admits an asymptotic average growth vector (also called here cycle time vector) $\chi$ 
  whose components are all equal to $\mu$.
\end{theorem}

\section{Max-plus algebra traffic model}
\label{sec-mp}

We present here a first traffic model for the train dynamics in a linear metro line.
The model describes the train dynamics without computing or applying any control of the train dwell times at platforms.
It simply writes the train departure times under two constraints on the inter-station travel time (running time + dwell time),
and on the safety time between successive trains. The model of this section does not take into account the effect of passengers on the
train dwell times at platforms.

\begin{figure}[thpb]
      \centering
	  \includegraphics[scale = 0.6]{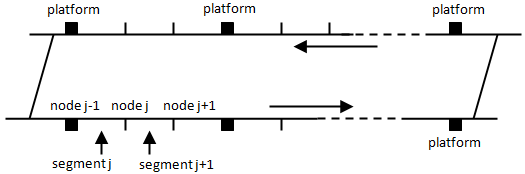}
      \caption{Representation of a linear metro line.}
      \label{fig-loop}
\end{figure}

Let us consider a linear metro line of $N$ platforms as shown in Figure~\ref{fig-loop}.
In order to model the train dynamics on the whole line, including the dynamics
on inter-stations, we discretize the inter-stations space, and thus the whole line, in segments
(or sections, or blocks). The length of every segment must be larger than the length of
one train. We then consider the following notations. \\~~

\noindent
\begin{tabular}{ll}
  $N$ & number of platforms.\\
  $n$ & number of all segments of the line.\\ 
  $m$ & number of running trains.\\
  $L$ & the length of the whole line. \\ 
  $b_j$ & $\in \{0,1\}$: boolean number of trains being on \\ 
        & segment $j$ at time zero.\\
  $\bar{b}_j$ & $= 1 - b_j \in \{0,1\}$.
\end{tabular}~~\\
\begin{tabular}{ll}  
  $d^k_j$ & instant of the $k$th departure from node $j$. \\
          & Notice that $k$ do not index trains, but count the \\
          & number of departures from segment $j$. \\
  $a^k_j$ & instant of the $k$th arrival to node $j$. \\
          & Notice that $k$ do not index trains, but count the \\
          & number of arrivals to segment $j$.
\end{tabular}~~\\
\begin{tabular}{ll}  	  
  $r_j$ & the running time of a train on segment $j$, i.e. from \\
        & node $j-1$ to node $j$.\\
  $w_j^k$ & $ = d_j^k - a_j^k$: train dwell time corresponding to \\
          & the $k$th arrival to- and departure from node $j$.
\end{tabular}~~\\
\begin{tabular}{ll}            
  $t_j^k$ & $ = r_{j} + w_j^k$: train travel time from node $j-1$ to\\
          & node $j$, corresponding to the $k$th arrival to- and \\
          & departure from node $j$.
\end{tabular}
\begin{tabular}{ll}  		
  $g^k_j$ & $ = a_{j}^{k} - d^{k-1}_{j}$: node- (or station-) safe separation\\
          & time (also known as close-in time), corresponding \\
  	  & to the $k$th arrival to- and $(k-1)$st departure\\
  	  & from node $j$.\\
  $h_j^k$ & $ = d^k_j - d^{k-1}_j = g_j^k + w_j^k$~: departure time headway at \\
          & node $j$, associated to the $(k-1)$st and $k$th \\
          & departures from node $j$.\\
  $s_j^k$ & $ = g_j^{k+b_j} - r_j^k$: a kind of node safe separation time\\
          & which does not take into account the running time.
\end{tabular}~~\\

We also use underlined and over-lined notations to denote the maximum and minimum bounds of 
the corresponding variables respectively.
Then $\bar{r}_j, \bar{t}_j, \bar{w}_j, \bar{g}, \bar{h}_j$ and $\bar{s}_j$ and respectively
$\underline{r}_j, \underline{t}_j, \underline{w}_j, \underline{g}, \underline{h}_j$ and $\underline{s}_j$
denote maximum and minimum running, travel, dwell, safe separation, headway and $s$ times, respectively.

The average on $j$ and on $k$ (asymptotic) of those variables are denoted without any subscript or superscript.
Then $r, t, w, g, h$ and $s$ denote the average running, travel, dwell, safe separation, headway and $s$ times,
respectively.

It is easy to check the following relationships.  
  \begin{align}
    & g = r + s, \label{form1} \\
    & t = r + w, \label{form2} \\
    & h = g + w = t + s = \frac{n}{m} t = \frac{n}{n-m} s. \label{form3}
  \end{align}    
Indeed, (\ref{form1}) comes from the definition of $s_j^k$ and~(\ref{form2}) comes from
the definition of $t_j^k$. For~(\ref{form3}),
\begin{itemize}
  \item $h=g+w$ comes from the definition of $h_j^k$,
  \item $h=t+s$ comes from the definition of $t_j^k$ and $s_j^k$ and from $h=g+w$,
  \item $h=nt/m$ average train time-headway is given by the travel time of the whole line ($nt$) divided by
     the number of trains.
  \item $h=ns/(n-m)$ can be derived from $h=t+s$ and $h=nt/m$.
\end{itemize}  

The running times of trains on inter-stations are assumed to be constant.
The running times $r_j$ of trains on every segment $j$, are also considered to be constant.
They can be calculated from the running times on inter-stations, and by means of given inter-station speed profiles,
depending on the characteristics of the line and of the trains running on it.
We then have $\underline{t}_j = r_j + \underline{w}_j$ and $\bar{t}_j = r_j + \bar{w}_j$ for every $j$.

An important remark here is that the variable $w_j^k$ denote dwell times at all nodes $j\in\{1,\ldots,n\}$
including non-platform nodes. The lower bounds $\underline{w}_j$ should be zero for the non-platform nodes $j$, and they
should be strictly positive for platform nodes. Therefore, the asymptotic average dwell time $w$ on all the nodes
is different from (it is underestimated comparing to) the asymptotic dwell time on the platform nodes.
We have the same remark also for the variables $g, t$ and $s$. In order to clarify this, we use the additional notations.
\begin{itemize}
  \item $w^*$: asymptotic average dwell time on platforms.
  \item $g^*$: asymptotic average safe separation time on platforms.
  \item $t^* = r + w^*$: asymptotic travel time on segment $j$ upstream of platform node $j$.
  \item $s^* = g^* - r$: asymptotic average safety time between node $j-1$ and platform node $j$.
\end{itemize}

Another important remark is that, as we consider constant running times on segments and on inter-stations,
then every train deceleration or stopping at the level of an inter-station, generally caused by an interaction with the train ahead,
is modeled here by a dwell time extension at one of the nodes at the considered inter-station.
We note that the inter-station train running times can also be considered as control variables, in addition to train dwell times at platforms.
In such modeling, shortened inter-station train running times can compensate extended train dwell times at platforms.
We shall consider this extended modeling in a future work. 

The model we present in this section is built on two time constraints:
\begin{itemize}
  \item A constraint on the travel time on every segment $j$.
     \begin{equation}\label{const1}
       d^k_j \geq d^{k-b_{j}}_{j-1} + \underline{t}_j.
     \end{equation}
     Constraint~(\ref{const1}) tells first that the $k$th departure from node $j$ corresponds 
     to the $k$th departure from node $(j-1)$ in case where there is no train at segment $j$
     at time zero ($b_j=0$), and corresponds to the $(k-1)$st departure from node $(j-1)$
     in case where there is a train at segment $j$ at time zero.
     Constraint~(\ref{const1}) tells in addition that the departure from node $j$
     cannot be realized before the corresponding departure from node $(j-1)$ plus
     the minimum travel time $\underline{t}_j = \underline{r}_j + \underline{w}_j$ from node $j-1$ to node $j$.
  \item A constraint on the safe separation time at every segment~$j$.
     $$\begin{array}{ll}
         d_j^k - d^{k-\bar{b}_{j+1}}_{j+1} & = a_{j+1}^{k+b_{j+1}} - r_{j+1} - d^{k-\bar{b}_{j+1}}_{j+1} \\~~\\
                                           & = g_{j+1}^{k+b_{j+1}} - r_{j+1} \\~~\\
                                           & \geq \underline{g}_{j+1} - r_{j+1} = \underline{s}_{j+1}.
       \end{array}$$
     That is
     \begin{equation}\label{const2}
        d^k_j \geq d^{k-\bar{b}_{j+1}}_{j+1} + \underline{s}_{j+1}.
     \end{equation}
     Constraint~(\ref{const2}) tells first that, in term of safety, the $k$th departure from node $j$
     is constrained by the $(k-1)$st departure from node $(j+1)$ in case where there is no train at segment
     $(j+1)$ at time zero, and it is constrained by the $k$th departure from node $(j+1)$ in case where there is
     a train at segment $(j+1)$ at time zero. Constraint~(\ref{const2}) tells in addition that the $k$th
     departure from node $j$ cannot be realized before the departure constraining it from node $(j+1)$
     plus the minimum safety time at node $(j+1)$.
\end{itemize}

The model then combines constraints~(\ref{const1}) and~(\ref{const2}), and gives the $k$th train departure time from each segment $j$, as follows.
\begin{equation}\label{eq-d1}
  d^k_j = \max \{d^{k-b_{j}}_{j-1} + \underline{t}_j, d^{k-\bar{b}_{j+1}}_{j+1} + \underline{s}_{j+1} \}, \; k\geq 1, 1\leq j\leq n,
\end{equation}
where the index $j$ is taken with modulo $n$. That is to say that, for the two particular cases of $j=1$ and $j=n$, the dynamics are written as follows.
\begin{align}
  & d^k_1 = \max \{d^{k-b_{1}}_n + \underline{t}_1, d^{k-\bar{b}_2}_{2} + \underline{s}_2 \}, \quad k\in \mathbb N, \nonumber \\
  & d^k_n = \max \{d^{k-b_{n}}_{n-1} + \underline{t}_n, d^{k-\bar{b}_1}_{1} + \underline{s}_1 \}, \quad k\in \mathbb N. \nonumber
\end{align}

The maximum operator of the dynamics~(\ref{eq-d1}) means that the $k^{\text{th}}$ departure time from node $j$ takes effect 
as soon as the two constraints~(\ref{const1}) and~(\ref{const2}) are satisfied. Equivalently, we can say that 
the $k^{\text{th}}$ departure time from node $j$ is given by the minimum instant satisfying the two constraints~(\ref{const1}) and~(\ref{const2}).

With Max-plus notations, and using the back-shift operator~$\gamma$, defined in 
section \ref{sec-review}, the dynamics~(\ref{eq-d1})
are written as follows.
\begin{equation}\label{eq-d11}
  d_j = \underline{t}_j \gamma^{b_{j}} d_{j-1} \oplus \underline{s}_{j+1} \gamma^{\bar{b}_{j+1}} d_{j+1}, \quad 1\leq j\leq n.
\end{equation}

We denote by $d^k$ the vector with components $d^k_j$ for $j=1,\ldots, n$. The dynamics~(\ref{eq-d11}) are then
written as follows.
\begin{equation}
   d^k = A(\gamma) \otimes d^k,
\end{equation}
where $A(\gamma)$ is the following Max-plus polynomial matrix.
\small
$$A(\gamma) = \begin{pmatrix}
                 \varepsilon & \gamma^{\bar{b}_2}\underline{s}_2 & \varepsilon & \cdots & \varepsilon & \gamma^{b_1}\underline{t}_1 \\
                 \gamma^{b_2}\underline{t}_2 & \varepsilon & \gamma^{\bar{b}_3}\underline{s}_3 & \varepsilon & \cdots & \varepsilon \\
                  & \ddots & \varepsilon & \ddots & & \\
                  \varepsilon & \cdots & \gamma^{b_{j}}\underline{t}_j & \varepsilon & \gamma^{\bar{b}_{j+1}}\underline{s}_{j+1} & \varepsilon \\
                   & & & \ddots & \varepsilon & \\
                  \gamma^{\bar{b}_1}\underline{s}_1 & \varepsilon & \cdots & \varepsilon & \gamma^{b_{n}}\underline{t}_n & \varepsilon
              \end{pmatrix}$$
\normalsize

\begin{theorem}\label{th-mpm}
  The train dynamics~(\ref{eq-d1}) converges to a stable stationary regime with a unique average asymptotic growth vector,
  whose components are all equal and are
  interpreted here as the asymptotic average train time-headway $h$ given as follows.
  $$h = \max \left\{ \frac{\sum_j \underline{t}_j}{m}, \max_j (\underline{t}_j+\underline{s}_j), \frac{\sum_j \underline{s}_j}{n-m}\right\}.$$
\end{theorem}~\\

\begin{figure}[thpb]
      \centering
	  \includegraphics[scale = 0.42]{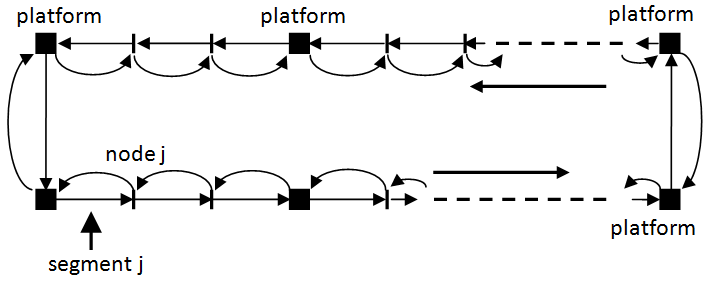}
      \caption{The graph $\mathcal G(A(\gamma))$.}
      \label{graph1}
\end{figure}

\proof
The graph $\mathcal G(A(\gamma))$ associated to the matrix $A(\gamma)$ is strongly 
connected; see Figure~\ref{graph1}.
Therefore, by Theorem~\ref{th-mpa}, we know that the asymptotic average growth 
vector of the dynamic system~(\ref{eq-d1}), whose components $h_j$
are interpreted here as the asymptotic average time-headway of the trains on segment $j$,
exists, and that all its components are the same
$h = h_j = \lim_{k\to +\infty} d^k_j/k, \forall j=1,\ldots,n$.
Moreover, $h$ coincides with the unique generalized eigenvalue of $A(\gamma)$,
given by Theorem~\ref{th-mpa} as the maximum cycle mean of the graph $\mathcal G(A(\gamma))$.
Three different elementary cycles are distinguished in $\mathcal 
G(A(\gamma))$; see Figure~\ref{graph1}.

\begin{itemize}
  \item The Hamiltonian cycle $c$ in the direction of the train running, with cycle mean
     $$\frac{W(c)}{D(c)} = \frac{\sum_j \underline{t}_j}{\sum_j b_j} = \frac{\sum_j \underline{t}_j}{m}.$$
  \item All the cycles $c_j$ of two links relying nodes $j-1$ and $j$, with cycle means
     $$\frac{W(c_j)}{D(c_j)} = \frac{(\underline{t}_j+\underline{s}_j)}{(b_j+\bar{b}_j)} = \underline{t}_j+\underline{s}_j, \quad \forall j.$$
  \item The Hamiltonian cycle $\bar{c}$ in the reverse direction of the train running, with cycle mean
     $$\frac{W(\bar{c})}{D(\bar{c})} = \frac{\sum_j \underline{s}_j}{\sum_j \bar{b}_j} = \frac{\sum_j \underline{s}_j}{n-m}.$$
\end{itemize}
\endproof

An important remark on Theorem~\ref{th-mpm} is that the asymptotic average 
train time-headway depends
on the average number of trains running on the metro line, without depending on 
the initial departure times of the trains (initial condition of the dynamic system).

\subsection{Fundamental traffic diagram for train dynamics}
\label{sec-ftd}

By similarity to the road traffic, one can define what is called \textit{fundamental traffic diagram} for the train dynamics.
In road traffic, such diagrams give relationships between car-flow and car-density on a road section; see for example~\cite{Far12, FHL13};
also extended to network (or macroscopic) fundamental diagrams; see for example~\cite{FGQ05, FGQ11, FGQ07, Far09, FGQ07b,FGQ11b}.

Let us first notice that the result given in Theorem~\ref{th-mpm} can be written as follows.
\begin{equation}\label{diag1}
  h(\sigma) = \max\left\{ \tau \sigma, h_{\min}, \frac{\omega}{\frac{1}{\underline{\sigma}} - \frac{1}{\sigma}}\right\},
\end{equation}
where $h$ is the asymptotic average train time-headway, $\sigma := L/m$ is the average train space-headway, $\tau := \sum_j \underline{t}_j / L = 1/v$
is the inverse of the maximum train speed $v$, $h_{\min} := \max_j h_j = \max_j (\underline{t}_j + \underline{s}_j)$,
$\omega := \sum_j \underline{s}_j / L$, and $\underline{\sigma} := L/n$ is the minimum train space-headway.
Relationship~(\ref{diag1}) gives the asymptotic average train time-headway as a function of the average train space-headway.

One can also write a relationship giving the average train time-headway as a function of the average train density $\rho := m/L = 1/\sigma$.
\begin{equation}\label{diag2}
  h(\rho) = \max\left\{ \frac{\tau}{\rho}, h_{\min}, \frac{\omega}{\bar{\rho} - \rho}\right\},
\end{equation}
where $\bar{\rho} := n/L = 1/\underline{\sigma}$ is the maximum train density on the metro line.

Let us know denote by $f=1/h$ the average train frequency (or flow) on the metro line, which is, as well known, the inverse of
the average train time-headway.
Then, from~(\ref{diag2}), we obtain a trapezoidal fundamental traffic diagram (well known in the road traffic) for the metro line.
\begin{equation}\label{diag3}
  f(\rho) = \min \left\{v \rho,  f_{\max}, w' (\bar{\rho} - \rho)\right\},
\end{equation}
where $f_{\max} = 1/h_{\min}$ is the maximum train frequency over the metro line segments,
$v = 1/\tau$ is the free (or maximum) train-speed on the metro line, and $w' = 1/\omega$ is the backward wave-speed
for the train dynamics.\footnote{We use the notation $w'$, with a prime, for the backward wave speed, in order to distinguish 
it with dwell time notation $w$.}

Relationships~(\ref{diag1}),~(\ref{diag2}) and~(\ref{diag3}) show how the asymptotic average train time-headway, 
and the asymptotic average train frequency change in function of the number of trains running on the metro line.
Moreover, they give the (maximum) train capacity of the metro line (expressed by the average train time-headway or by the average train frequency),
as well as the corresponding optimal number of trains. Furthermore, those relationships describe wholly the traffic phases of the train dynamics.

\begin{theorem}\label{cor-wg}
  The asymptotic average dwell time $w$ and safe separation time $g$ are given as follows.
  \begin{align}
    & w(\rho) = \max\left\{\underline{w}, \frac{h_{\min}}{\bar{\rho}}\;\rho - r, \frac{\omega}{\bar{\rho} - \rho} - \underline{g} \right\}.\label{diag-w}\\
    & g(\rho) = \max\left\{ \frac{\tau}{\rho} - \underline{w}, (r+h_{\min}) - \frac{h_{\min}}{\bar{\rho}}\;\rho, \underline{g} \right\}. \label{diag-g}
  \end{align}  
  where $\underline{w} = \sum_j \underline{w}_j /n, r = \sum_j r_j /n$ and $\underline{g} = \sum_j \underline{g_j} /n$.
\end{theorem}
\proof
We have
\begin{itemize}
  \item By (\ref{form2}) and (\ref{form3}), $w = t - r = (m/n) h - r$, then we replace $h$ using~(\ref{diag2}).
  \item By (\ref{form1}) and (\ref{form3}), $g = r + s = r + ((n-m)/n) h$, then we replace $h$ using~(\ref{diag2}). 
    Or directly form~(\ref{form3}), we have $g = h - w$, then we replace $h$ using~(\ref{diag2}) and $w$ using~(\ref{diag-w}). \endproof
\end{itemize}

\begin{figure}[thpb]
      \begin{center}    
      \includegraphics[scale=0.35]{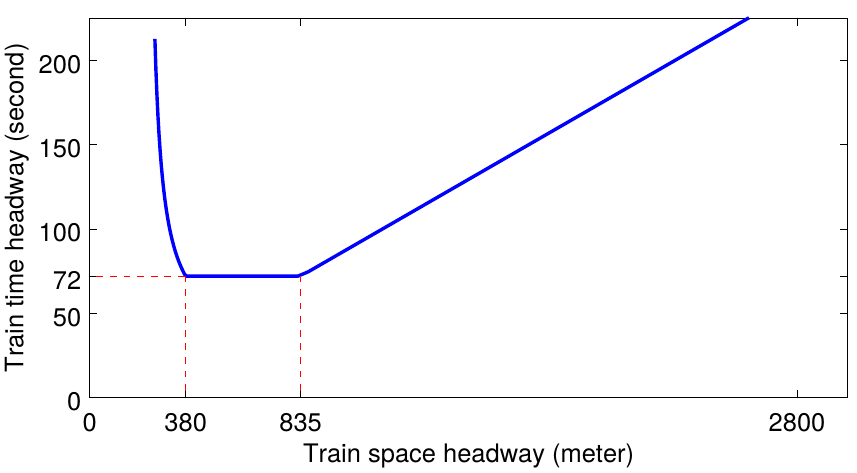} \\ ~~ \\
      \includegraphics[scale=0.35]{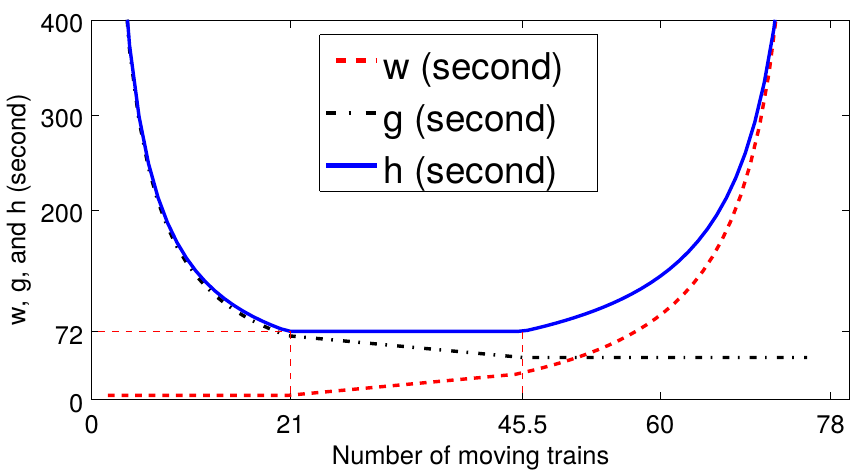} \\ ~~ \\
      \includegraphics[scale=0.35]{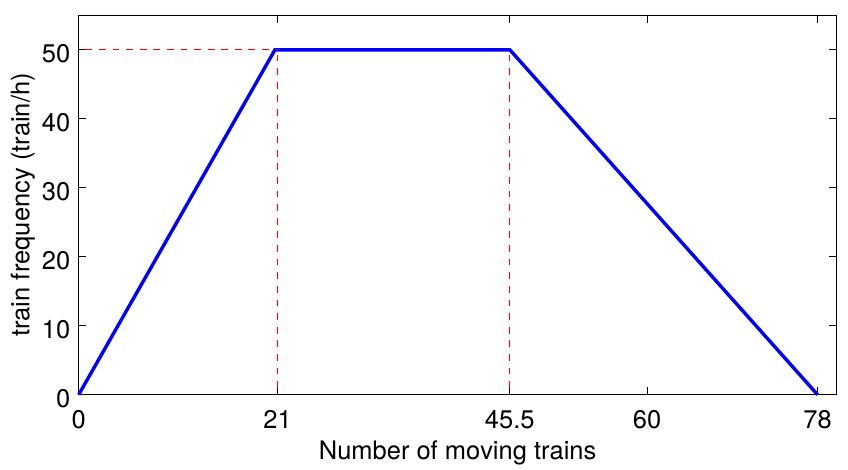}
      \caption{Analytical Phase diagrams for the train dynamics in a linear metro line.}
      \label{fig-diag1}
      \end{center}
\end{figure}

Figure \ref{fig-diag1} illustrates the relationships~(\ref{diag1}),
(\ref{diag2}), (\ref{diag3}), (\ref{diag-w}) and (\ref{diag-g}) for a linear metro line of 9 stations (18 platforms),
inspired from the automated metro line~14 of Paris~\cite{ratp99}.
The parameters considered for the line are given in Table~\ref{tab-param}.

\begin{table}
\centering
\caption{Parameters of the metro line considered.}
\begin{tabular}{|l||l|}
  \hline
  Number of stations & 9 ($\Rightarrow $ 18 platforms)\\
  \hline
  Segment length & about 200 meters (m.) \\
  \hline
  Free train speed $v_{\text{run}}$ & 22 m/s (about 80 km/h) \\
  \hline 
  Train speed on terminus & 11 m/s (about 40 km/h) \\
  \hline
  Min. dwell time $\underline{w}$ & 20 seconds \\
  \hline
  Min. safety time $\underline{s}$ & 30 seconds \\
  \hline  
  Inter-station length (in meters) & \begin{tabular}{lr}
			       S.-Laz. $\to$ Mad. & 618 m. \\
			       Mad. $\to$ Pyr. & 712 m. \\
                               Pyr. $\to$ Cha. & 1359 m. \\
			       Cha. $\to$ G.-Lyo & 2499 m.\\
			       G.-Lyo $\to$ Ber. & 624 m.\\
			       Ber. $\to$ C. S. Emi. & 970 m. \\
			       C. S. Emi. $\to$ Bib. & 947 m.\\
			       Bib. $\to$ Oly. & 713 m.
			     \end{tabular} \\
  \hline                          
\end{tabular}
\label{tab-param}
\end{table}

According to Figure~\ref{fig-diag1}, the maximum average train frequency for the considered metro line is about 50 trains/hour,
corresponding to an average time-headway of 72 seconds. The optimal number of running trains to reach the
maximum frequency is 21 trains. We note that time-margins for robustness are not considered here. 

Formulas~(\ref{diag-w}) and~(\ref{diag-g}) are important for the control model we present 
in section~\ref{sec-stable}, where we consider $w$ as the control vector and $g$ as the
traffic state vector, and where the formulas~(\ref{diag2}),~(\ref{diag-w}), and~(\ref{diag-g})
of the max-plus traffic model are used.

\subsection{The traffic phases of the train dynamics}
\label{sec-tph}

Theorems~\ref{th-mpm} and~\ref{cor-wg}, and formulas~(\ref{diag1}), (\ref{diag2}) and~(\ref{diag3}) allow the description of the
traffic phases of the train dynamics~\ref{eq-d1}. Three traffic phases are distinguished.

\textit{\textbf{Free flow traffic phase}}. ($0 \leq \rho \leq f_{\max}/v$). During this phase, trains move freely on the line, which
    operates under capacity, with high average train time-headways. The average time-headway is a
    sum of the average minimum train dwell time with the average train safe separation time.
    The average train dwell time is independent of the number of
    running trains, while the average train time-headway as well as the average train safe separation time
    decrease rapidly with the number of running trains. We notice that the average train frequency increases linearly with respect to 
    the number of running trains. Similarly, the average train time-headway increases linearly with respect to the average space-headway.

\textit{\textbf{Maximum train-frequency traffic phase}}. ($f_{\max}/v \leq \rho \leq \bar{\rho} - f_{\max}/w'$).
    During this phase, the metro line operates at its maximum train-capacity. The latter is constant, i.e. independent of the
    number of running trains. The average train dwell time $w$ increases linearly with the number of the running trains.
    The average train safe separation time $g$ decreases linearly with the number of running trains. 
    The average train time-headway $h = g + w$ remains constant and independent of the number of running trains.
    The optimum number of running trains on the line is attained at the beginning of the this traffic phase.
    This optimum number is $L f_{\max}/v$, corresponding to train density $f_{\max}/v$.
    However, in the case where passenger arrivals are taken into account, it can be interesting to increase the number of running
    trains on the line. Indeed, although the average train time-headway remains constant during this phase, the average
    train dwell time increases, while the average train safe separation time decreases with the increasing of the number of trains
    running on the line. This induces less average train safe separation time, so less time for the accumulation of passengers on platforms
    in one side, and more time to passengers to go onto the trains on the other side,
    without affecting the average train time-headway; see Figure~\ref{fig-diag1}.

\textit{\textbf{Congestion traffic phase}}. ($\bar{\rho} - f_{\max}/w' \leq \rho \leq \bar{\rho}$).
    During this phase, trains interact with each other and the metro line operates under capacity with high average train time-headways.
    The average train time-headway is given by the sum of the average train dwell time with the average minimum safe separation time.
    The latter is independent of the number of running trains, while the average train time-headways, as well as the average train dwell times
    increase rapidly with the number of running trains on the metro line. We notice that the average train frequency
    decreases linearly with the number of trains running on the metro line.

Before introducing passenger arrivals in the train dynamics, we give here an idea on their effect in term of the service offered.
Let us consider the metro line as a server of passengers, with an average passenger arrival rate $\lambda$ by platform.
Under the assumption of unlimited passenger capacity of trains, the average service rate of passengers by platform can be fixed to $\alpha w^* /h$, where
$\alpha$ is the average passenger upload rate by platform, $w^*$ denotes, as indicated above, the average dwell time on a platform,
and $h$ is the average train time-headway of the metro line.
In this case, the system, as a server of passengers, is stable under the following condition.
\begin{equation}\label{stab-cond}
   \lambda < \alpha w^* /h.
\end{equation}   

\section{Traffic instability due to passengers}
\label{sec-tip}

We consider in this section, train dynamics that take into account the passengers' travel demand.
We introduce this dependence by assuming that the train dwell time $w_j$ at platform $j$ depends
on the passenger volume at platform $j$, which depends on the safe separation time $g_j$ on the same platform.
In this section, we present a first model that extends naturally the model of section~\ref{sec-mp},
without optimizing the effect of the passengers volumes on the train dynamics.
In other words, we assume that the dependence of the dwell times at platforms with the passenger
demand is not controlled. We will see that in this case, the dynamic system is unstable.

We do not consider in this article a dynamic model for the number of passengers at platforms.
Therefore, we assume that the train dwell times at platforms depend directly on the passenger arrival rates.
We model here the effect of the embarking time of passengers into the trains, on the train dwell time.
We do not model the effect of the alighting time of passengers from the trains, on the train dwell time.

In order to model the effect of passengers on the train dwell times at platforms, we consider the following
additional constraint on the dwell time at platforms.
\begin{equation}\label{eq-w1}
  w_j^k \geq \begin{cases}
                \frac{\lambda_j}{\alpha_j} g_j^k, & \text{ if } j \text{ indexes a platform}, \\
                0 & \text{ otherwise}.
             \end{cases}
\end{equation}
where $\alpha_j$ is the total passenger upload rate from platform $j$ onto the trains, if $j$ indexes a platform;
and $\alpha_j$ is zero otherwise;
and $\lambda_j$ is the average rate of the total arrival flow of passengers to platform $j$, if $j$ indexes a platform;
and $\lambda_j$ is zero otherwise.
$$\lambda_j = \left\{ \begin{array}{ll}
		                 \sum_i \lambda_{ji} & \text{ if } j \text{ indexes a platform} \\
		                 0 & \text{otherwise}.
		              \end{array} \right.$$
$\lambda_{ji}$ denote here the origin-destination arrival rates of passengers to platform $j$, with platform $i$ as destination.

By taking into account the additional constraint~(\ref{eq-w1}), the constraint~(\ref{const1}) is modified as follows.
\begin{equation}\label{new-constr}
	  d_j^k \geq d_{j-1}^{k-b_{j}} + r_j +  \max \left\{ \underline{w}_j, \frac{\lambda_j}{\alpha_j} g_j^k \right\}.
\end{equation}

We now modify the dynamic model~(\ref{eq-d1}) by taking into account the new constraint~(\ref{new-constr}).
We obtain, for nodes $j$ indexing platforms, the following.

\begin{equation}\label{eq-dd}
  d_j^k = \max \left\{
  \begin{array}{l}
     d_{j-1}^{k-b_{j}} + r_{j} + \underline{w}_j, \\ ~~\\
     \left(1+\frac{\lambda_j}{\alpha_j}\right) d_{j-1}^{k-b_{j}} - \left(\frac{\lambda_j}{\alpha_j}\right) d_j^{k-1} + \left(1+\frac{\lambda_j}{\alpha_j}\right) r_{j}, \\~~\\
     d_{j+1}^{k-\bar{b}_{j+1}} + \underline{s}_{j+1}.
  \end{array} \right.
\end{equation}
For non-platform nodes, the dynamics remain as in~(\ref{eq-d1}).

Let us notice that the dynamic system~(\ref{eq-dd}) has explicit and implicit terms.
Moreover, it can be written in the form~(\ref{eq-dp3}), i.e. as follows.
\begin{equation}\label{matrix_form}
  d_j^k = \max_{u\in\mathcal U} [(M^u d^{k-1})_j + (N^u d^{k})_j + c^u_j],
\end{equation}
where $M^u$ and $N^u$ are square matrices, and $c^u$ is a family of vectors, for $u\in \mathcal U$.
The matrices $N^u, u\in \mathcal U$ express implicit terms.
We notice that it suffices that $\exists j, \lambda_j > 0$ to have one of the matrices $M^u, u\in\mathcal U$ or $N^u, u\in\mathcal U$ not being 
sub-stochastic~\footnote{A matrix $A$ is sub-stochastic if $0 \leq A_{ij}\leq 1, \forall i,j$ and if $\sum_i A_{ij} \leq 1, \forall i$.}, 
since we have in this case $\alpha_j/(\alpha_j-\lambda_j) > 1$ and $-\lambda_j/(\alpha_j-\lambda_j) < 0$; see the dynamics~(\ref{eq-dd}).
Therefore, the dynamic system~(\ref{eq-dd}) cannot be seen as a dynamic programming system of a stochastic optimal control problem
of a Markov chain; see section~\ref{sub-sec-dps}.

It is easy to see that if $m=0$ or $m=n$, then the dynamic system~(\ref{eq-dd}) is fully implicit (it is not triangular), and it admits an
asymptotic regime with a unique additive eigenvalue $h = +\infty$, which is also the average growth rate of the system, and which is the asymptotic average 
train time-headway.
This case corresponds to $0$ or $n$ trains on the metro line. No train departure is possible for these two cases. We have the
average train frequency $f=0$ corresponding to the average time headway $h=+\infty$ (the additive eigenvalue of the system). 
One can also show that if $0 < m < n$, then the dynamic system~(\ref{eq-dd}) is triangular. That is to say that it is not fully implicit, and
there exists an order of updating the components of the state vector $d^k$, in such a way that no implicit term appears.

We know that the dynamic system~(\ref{eq-dd}) is not stable in general.
We justify this assertion with three arguments, two theoretical ones, and an application interpretation one.

First, as we noticed above, the matrices $M^u, u\in\mathcal U$ or $N^u, u\in\mathcal U$ are not sub-stochastic.
   The consequence of this is that the dynamics are not non-expansive. In fact, they are additive 1-homogeneous but not monotone.
   Many behaviors are possible for the dynamic system, depending on the parameters and on the initial
   conditions (expansive behavior, chaotic behavior, etc.)

Second, let us consider the metro line as a server of passengers, with an average passenger arrival
   rate $\lambda$ to every platform. 
   The average service rate of passengers by platform can be fixed to $\alpha w^* /h$ (under the assumption of
   infinite passenger capacity of trains), where, as indicated above, 
   $\alpha$ is the average passenger upload rate by platform, $w^*$ is the average dwell time at a platform,
   and $h$ is the average train time-headway of the metro line.
   The system, as a server, is then stable under the condition~(\ref{stab-cond}).
   In the model considered in this section, and in the passenger congestion case, where
   the arrival rates are sufficiently high so that the second term of the maximum operator of the
   dynamics~(\ref{eq-dd}) is activated, we get $w^* = (\lambda/\alpha) g$. 
   Therefore, $\lambda = \alpha w^* /g > \alpha w^* /h$ since $g < h$.   
   That is to say that the average arrival rate exceeds the average service one; which is contradictory with~(\ref{stab-cond}).
   The metro line as a server is then unstable.
   
Third, in practice, let us assume the $k^{\text{th}}$ arrival of a given train to platform $j$ is delayed.   
   Then $g^k_j$ will increase by definition, and $w^k_j$ will also increase by application of~(\ref{eq-w1}).
   By consequent of the increasing of $w^k_j$, the departure $d^k_j$ will be delayed, and the delay
   of $d^k_j$ would be longer than the one of $a^k_j$, because, it accumulates the delay of $a^k_j$ 
   and the increasing of $w^k_j$. Consequently, the arrival of the same train to platform $j+1$
   (downstream of platform $j$) will be delayed longer comparing to its arrival to platform~$j$.
   Therefore, the application of the control law~(\ref{eq-w1}) amplifies the train delays and propagates
   them through the metro line. 

The instability of this kind of dynamics has already been pointed out; see for example~\cite{BCB91}.

\section{Stable dynamic programming model}
\label{sec-stable}

\subsection{The model}

In this section we modify the train dynamics~(\ref{eq-dd}) in order to guarantee its stability.
We know that system~(\ref{eq-dd}) is unstable because of the relationship~(\ref{eq-w1}).
In order to deal with this instability, we propose to replace the dwell time control formula~(\ref{eq-w1}), by the following.
\begin{equation}\label{eq-w4}
  w_j^k \geq \begin{cases}
                \overline{w}_j - \frac{\theta_j^k}{\lambda_j^k/\alpha_j^k} g_j^k & \text{ if } j \text{ indexes a platform}, \\
                0 & \text{ otherwise}.
             \end{cases} 
\end{equation}
where we reversed the sign of the relationship between the dwell time $w_j^k$ and the safe separation time $g_j^k$,
without reversing the relationship between the dwell time $w_j^k$ and the ratio $\lambda_j^k/\alpha_j^k$; and where
$\overline{w}_j$ (maximum dwell time on node $j$) and $\theta_j^k$ are control parameters to be fixed.

The dynamics~(\ref{eq-dd}) are now rewritten, for nodes $j$ indexing platforms, as follows.
\begin{equation}\label{eq-dd2}
  d_j^k = \max \left\{
  \begin{array}{l}
     d_{j-1}^{k-b_{j}} + r_{j} + \underline{w}_j, \\~~ \\
     \left(1 - \delta_j^k\right) d_{j-1}^{k-b_{j}} + \delta_j^k d_j^{k-1} + \left(1 - \delta_j^k\right) r_{j} + \overline{w}_j, \\ ~~\\
     d_{j+1}^{k-\bar{b}_{j+1}} + \underline{s}_{j+1},
  \end{array} \right.
\end{equation}
where $\delta_j^k = \theta_j^k \alpha_j^k/\lambda_j^k, \forall j, k$.
For non-platform nodes, the dynamics remain as in~(\ref{eq-d1}).

If $\delta_j^k$ are independent of $k$ for every $j$, then
the dynamic system (\ref{eq-dd2}) can be written in the form~(\ref{matrix_form}).

As for the dynamic system~(\ref{eq-dd}), if $m=0$ or $m=n$, then the system~(\ref{eq-dd2}) is fully implicit, and admits an
asymptotic regime with a unique average asymptotic train time-headway $h = +\infty$ (no train movement is possible).
If $0 < m < n$, then the dynamic system~(\ref{eq-dd2}) is triangular.
In this case, and if $\delta_j^k$ are independent of $k$ for every $j$, then the explicit form of the dynamic
system~(\ref{eq-dd2}) can be written in the form~(\ref{eq-dp3}).
If, in addition, $0 \leq \delta_j \leq 1, \forall j$, then $M^u$ and $N^u$ are sub-stochastic matrices,
i.e. satisfying $M^u_{ij} \geq 0, N^u_{ij} \geq 0$, $\sum_j (M^u_{ij}) \leq 1$ and $\sum_j (N^u_{ij}) \leq 1$.
Moreover, we have $\sum_j (M^u_{ij}+N^u_{ij}) = 1, \forall i,u$.
In this case,~(\ref{eq-dd2}) is a dynamic programming system of an optimal control problem of a Markov chain,
whose transition matrices
can be calculated from $M^u, u\in\mathcal U$ and $N^u, u\in\mathcal U$ (they are the matrices corresponding to the
equivalent explicit dynamic system
obtained by solving the implicit terms of~(\ref{eq-dd2})); and whose reward vectors are $c^u, u\in\mathcal U$.

\begin{theorem}\label{stable}
  If $0 < m < n$ and if $\delta_j^k$ are independent of $k$ for every $j$, and $0 \leq \delta_j \leq 1, \forall j$,
  then the dynamic system~(\ref{eq-dd2}) 
  admits a stationary regime with a unique additive eigenvalue $h$.
  Moreover, the system admits a unique asymptotic average growth vector whose components are all equal to $h$
  (interpreted here as the asymptotic average train time-headway), independent of the initial state $d^0$.
\end{theorem}

\begin{figure}[htpb]
      \centering
	  \includegraphics[scale = 0.42]{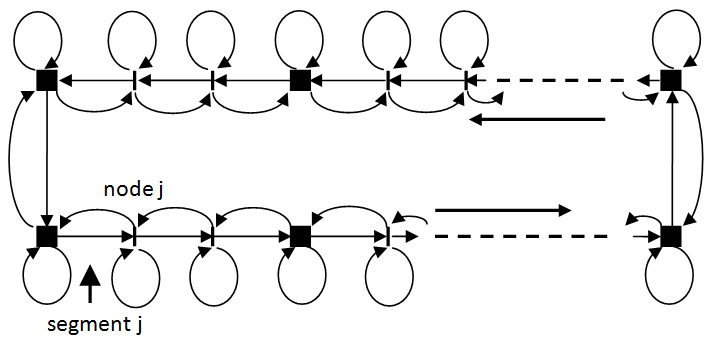}
      \caption{The graph associated to $\mathbf{f}$ in the proof of Theorem~\ref{stable}, which is also the same
         graph associated to the Max-plus linear system in the proof of Theorem~\ref{th-dp-mp}.}
      \label{graph2}
\end{figure}

\proof The proof uses Theorem~\ref{th-dps2}. 
  Let us denote by $\mathbf{f}$ the map associated to the dynamic system~(\ref{eq-dd2}), and use the notations
  $\mathbf{f}_i, i\in{0,1}$ and $\tilde{\mathbf{f}}$ as defined in section~\ref{sub-sec-gdps}.
  Since $M^u$ and $N^u$ are sub-stochastic matrices, with $\sum_j (M^u_{ij}+N^u_{ij}) = 1, \forall i,u$, then
  $\mathbf{f}$ is additive 1-homogeneous and monotone.
  The graph associated to $\mathbf{f}$ is strongly connected; see Figure~\ref{graph2}.
  Indeed, this graph includes the one of the Max-plus linear dynamics (graph of Figure~\ref{graph1}), which is already 
  strongly connected. The graph associated to $\mathbf{f}_0$ is acyclic since $0<m<n$.
\endproof  

We do not yet have an analytic formula for the asymptotic train time-headway (which we know coincides with eigenvalue $h$),
but Theorem~\ref{stable} guarantees its existence and its uniqueness.
Therefore, by iterating the dynamics~(\ref{eq-dd2}), one can approximate, for any fixed train density $\rho$, the associated asymptotic average
train time-headway $h(\rho)$ as follows.
\begin{equation}\label{hsim}
  h(\rho) \approx d^K_j/K, \forall j, \text{ for a large } K.
\end{equation}

In the following, we give a result (Theorem~\ref{th-dp-mp}) which tells us under which conditions on the control parameters
$\bar{w}_j$ and $\theta_j$ (or equivalently $\delta_j$), the dynamic system~(\ref{eq-dd2}) is a Max-plus linear system.
We will use this result in section~\ref{sub-sec-param}, in order to derive a general approach for fixing the control parameters
$\bar{w}_j$ and $\theta_j$ in such a way that the effect of passenger arrivals on the train dynamics will be well modeled.

\begin{theorem}\label{th-dp-mp}
  Let $\tilde{h}$ be the asymptotic average growth rate of the Max-plus linear system~(\ref{eq-d1}).
  The dynamic programming system~(\ref{eq-dd2}) with parameters $\bar{w}_j = \tilde{h}, \forall j$, and $\delta_j = 1, \forall j$,
  is a Max-plus linear system, whose asymptotic average growth rate coincides with $\tilde{h}$.
\end{theorem}

\proof It is easy to see that if $\delta_j = 1, \forall j$, then system~(\ref{eq-dd2}) is a Max-plus linear
system whose associated graph has $n$ additional cycles (which are loop-cycles) comparing to the one associated
to system~(\ref{eq-d1}); see Figure~\ref{graph2}. Moreover, if $\bar{w}_j = \tilde{h}, \forall j$, then
the cycle mean of the loops are all equal to $\tilde{h}$.
All the other parameters corresponding to the characteristics of the metro line and of the
trains running on it, remain the same as the ones of system~(\ref{eq-d1}).
Therefore, as in Theorem~\ref{th-mpm}, the asymptotic average  time-headway $h$ is given by the maximum cycle mean
of the graph associated to the Max-plus linear system obtained from system~(\ref{eq-dd2}).
Four different elementary cycles are distinguished on that graph; see Figure~\ref{graph2}.
\begin{itemize}
  \item The Hamiltonian cycle in the direction of the train movements, with mean $\sum_j \underline{t}_j/m$.
  \item All the cycles of two links relying nodes $j-1$ and $j$, with mean $\underline{t}_j+\underline{s}_j$ each.
  \item The Hamiltonian cycle in reverse direction of the train dynamics, with mean
    $\sum_j \underline{s}_j/(n-m)$.
  \item The $n$ loop-cycles with mean $\tilde{h}$.
\end{itemize}
Hence
$$h = \max \left\{ \frac{\sum_j \underline{t}_j}{m}, \max_j (\underline{t}_j+\underline{s}_j),
   \frac{\sum_j \underline{s}_j}{n-m}, \tilde{h}\right\} = \tilde{h}.$$
\endproof

\subsection{How to fix the control parameters $\bar{w}_j$ and $\theta_j$}
\label{sub-sec-param}

Theorem~\ref{th-dp-mp} tells us that if we fix $(\bar{w}_j, \theta^k_j) = (\tilde{h}(\rho), \lambda^k_j/\alpha^k_j)$
in~(\ref{eq-dd2}), or equivalently $(\bar{w}_j, \delta^k_j) = (\tilde{h}(\rho),1)$,
then, we obtain a Max-plus linear dynamic system, which does not take into account the passenger demand.

We consider here that $\alpha_j^k, \lambda_j^k, \theta_j^k$ and then $\delta_j^k$ are independent of $k$,
and then denote $\alpha_j, \lambda_j, \theta_j$ and $\delta_j$ respectively.
Let us consider as above the metro line as a server of passengers. We assume that the average service rate
at the level of one platform is $\alpha w^* /h$.
Under the Max-plus linear model for the train dynamics, the asymptotic average service rate depends on the number $m$ of trains,
or equivalently on the train density $\rho$. It is given by $\alpha w^* /\tilde{h}(\rho)$.
The average arrival of passengers to the platforms is $\lambda$.
Therefore, the metro line as a server operating under the Max-plus linear model is stable under the condition 
$\lambda < \alpha w^* /\tilde{h}(\rho)$.
Let use the following notation.
\begin{equation}
  \tilde{\lambda}_j(\rho) := \alpha w^* /\tilde{h}(\rho).
\end{equation}
The stability condition of the server means that if $\lambda_j \leq \tilde{\lambda}_j(\rho), \forall j$, then the server
operating under the Max-plus model
can serve the passengers without adapting the train dwell times to the passenger arrival rate.
Basing on this remark, we propose here an approach of fixing the parameters $\bar{w}$ and $\theta$ as functions of the average
arrival rate of passengers, in such a way that
\begin{itemize}
  \item If $\lambda_j \leq \tilde{\lambda}_j(\rho), \forall j$, then the dynamic system behaves as a Max-plus linear one.
    That is, the train dwell times are not constrained by the arrival rates of passengers.
  \item If $\exists j, \lambda_j > \tilde{\lambda}_j(\rho)$ then the system switches to a dynamic programming one, where
    the train dwell times are constrained by the arrival rates of passengers.
\end{itemize}

In order to have that, we fix the parameters $\bar{w}$ and $\theta$ as follows.
\begin{align}
  & \bar{w}_j(\rho) := \tilde{h}(\rho), \forall \rho, j. \label{param1} \\
  & \theta_j(\rho) := \frac{\tilde{\lambda}_j(\rho)}{\max\left(\lambda_j,\tilde{\lambda}_j(\rho)\right)} \; \frac{\lambda_j}{\alpha_j}, \forall \rho, j. \label{param2}
\end{align}
Let us notice that fixing $\theta(\rho)$ as in~(\ref{param2}) is equivalent to fixing $\delta(\rho)$ as follows.
\begin{equation}
   \delta_j(\rho) := \frac{\tilde{\lambda}_j(\rho)}{\max\left(\lambda_j,\tilde{\lambda}_j(\rho)\right)}, \forall \rho, j. \label{param3}
\end{equation}

We then have $0\leq \delta_j(\rho) \leq 1$ by definition, and
\begin{itemize}
  \item If $\lambda_j \leq \tilde{\lambda}_j(\rho), \forall j$, then $\delta_j(\rho) = 1, \forall j$, and the dynamic system behaves as a Max-plus linear one,
     and the train dwell times are not constrained by the arrival rates of passengers.
  \item If $\exists j, \lambda_j > \tilde{\lambda}_j(\rho)$ then $\exists j, \delta_j(\rho) < 1$, and the system switches to a dynamic programming one, and
     the train dwell times are constrained by the arrival rates of passengers.
\end{itemize}

\begin{table*}[thbp]
\centering
\caption{Asymptotic average train time-headway and frequency as functions of the number of running trains.
         For the symmetric passenger arrivals, the passenger arrival rate to any platform is equal to 1 times a factor $c$ given in the figure.
         For the asymmetric passenger arrivals, the average passenger arrival rates to platforms are given in Figure~\ref{arrival-rates} 
         times a factor $c$ given in the figure.}
\begin{tabular}{|c|c|}
  \hline
  Symmetric passenger arrival & Asymmetric passenger arrival \\
  \hline
   & \\
  \includegraphics[scale=0.16]{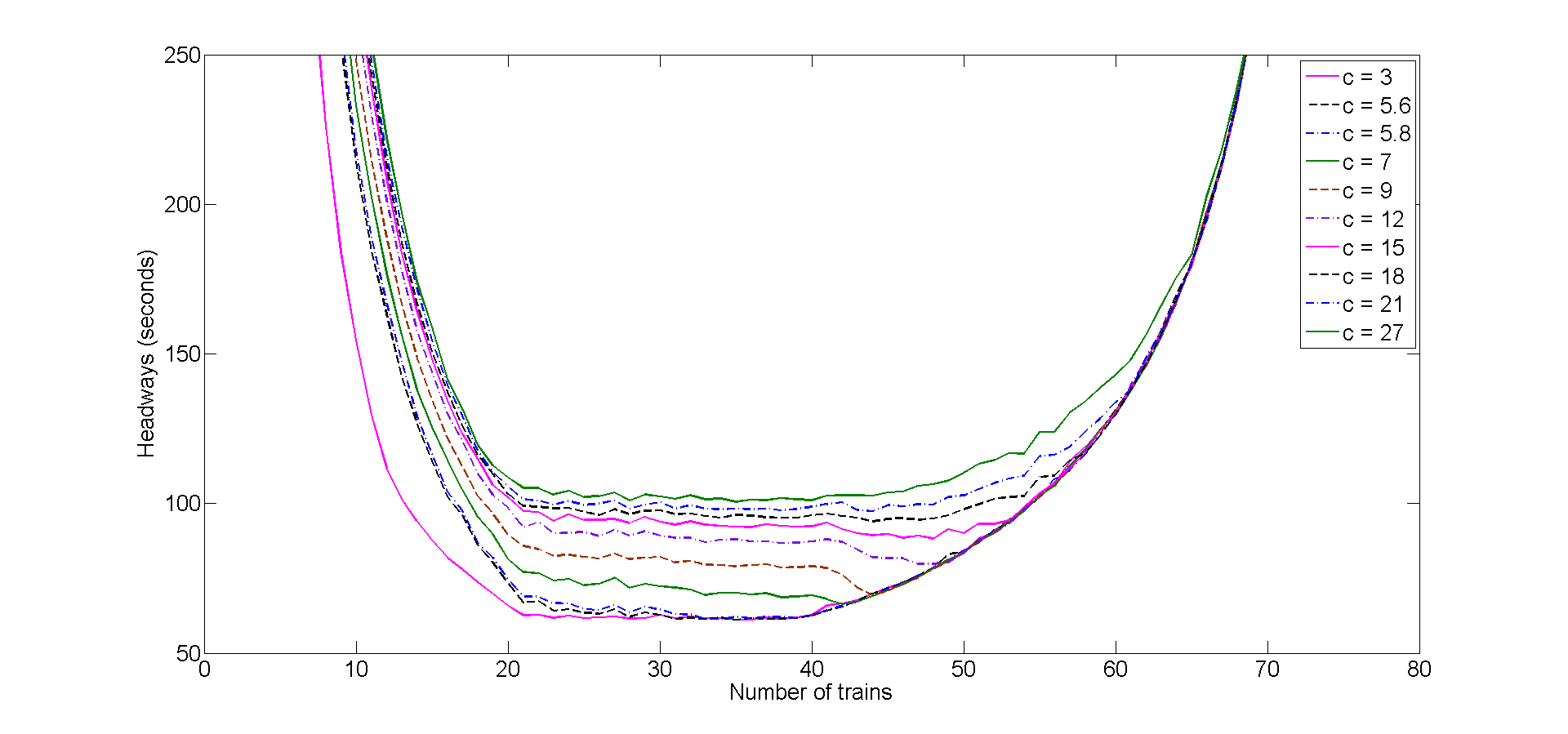} & \includegraphics[scale=0.16]{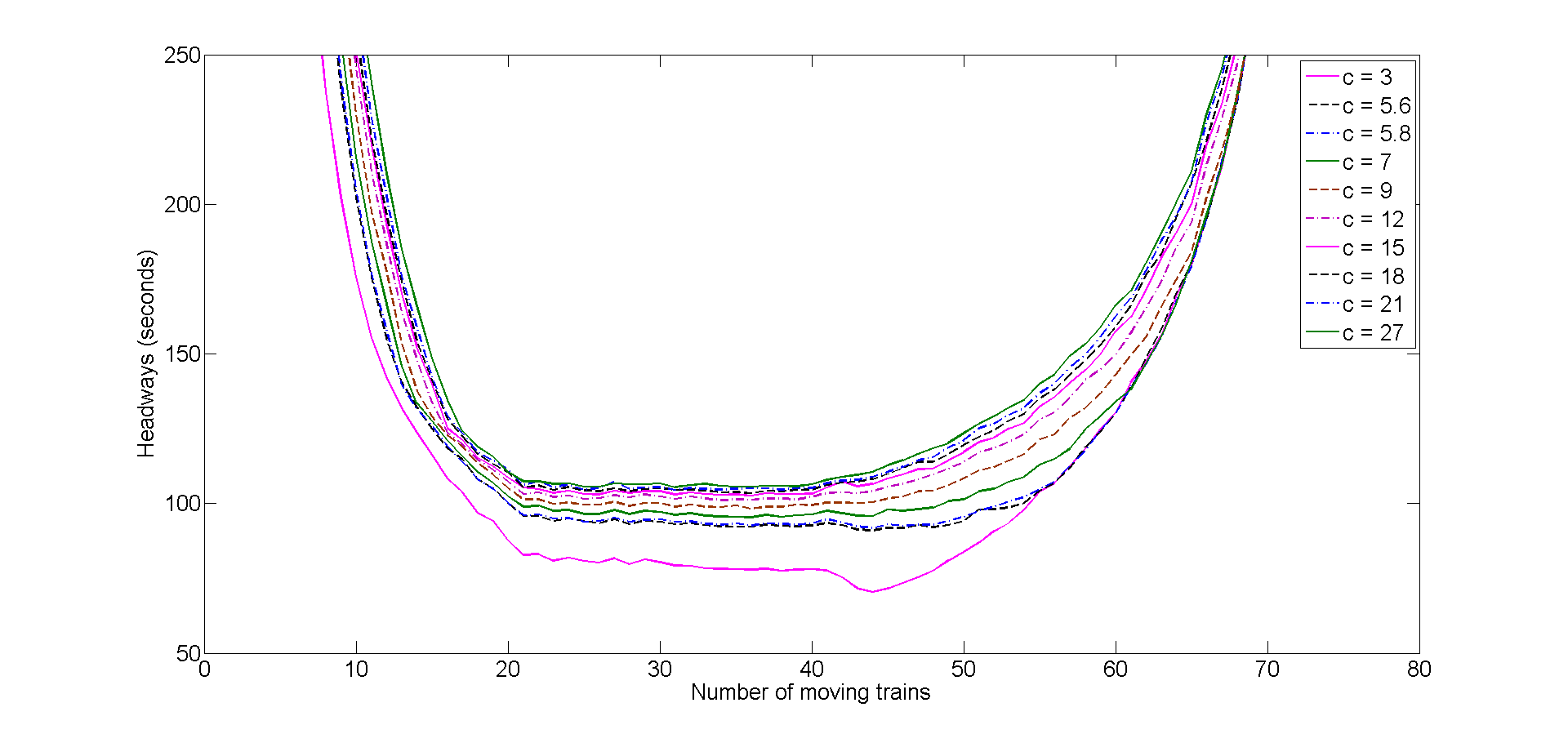} \\
  \includegraphics[scale=0.16]{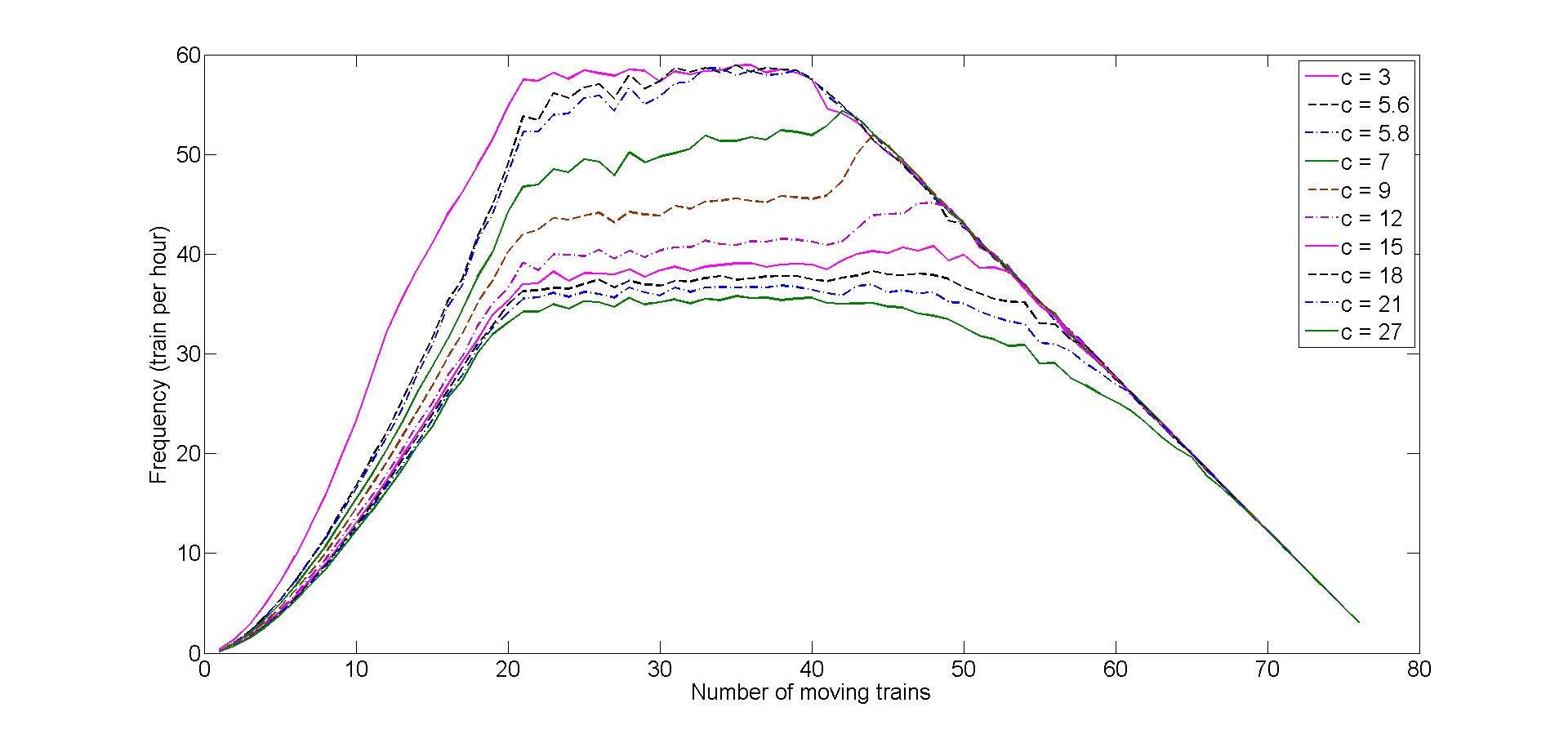} & \includegraphics[scale=0.16]{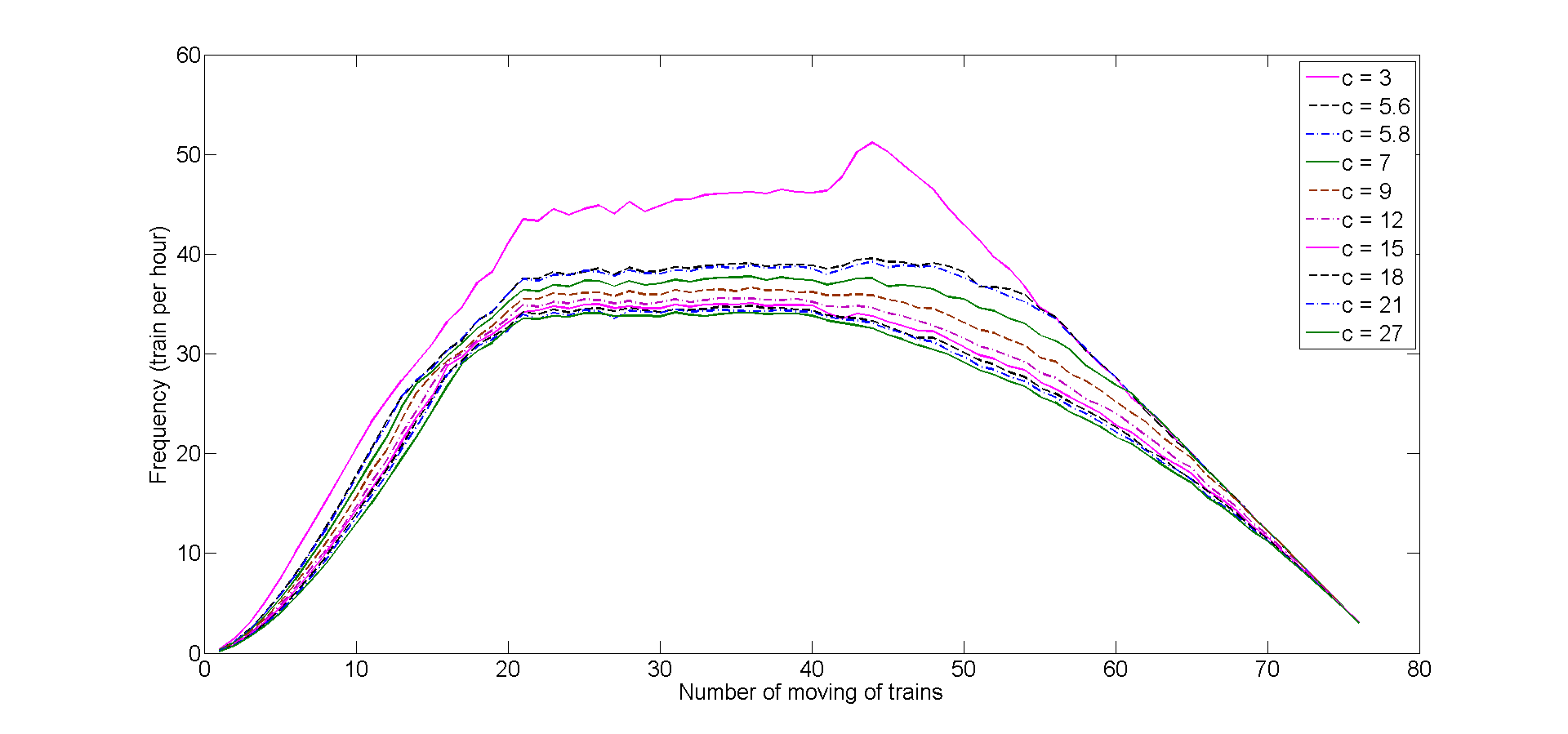} \\
  \hline
\end{tabular}
\label{tab_1}
\end{table*}

We summarize the latter findings in the following result (Theorem~\ref{th_stab}). \\~~
\begin{theorem}\label{th_stab}
  For any fixed value of the train density $\rho$ on the metro line,
  the dynamic system~(\ref{eq-dd2}) with parameters $\bar{w}_j$ and $\delta_j$ fixed dependent on $\rho$ as in~(\ref{param1}) and~(\ref{param3})
  respectively, 
  admits a stationary regime with a unique additive eigenvalue $h(\rho)$.
  Moreover, the system admits a unique asymptotic average growth vector whose components are all equal to $h(\rho)$
  (interpreted here as the asymptotic average train time-headway), independent of the initial state $d^0$.
  We have $h(\rho) \geq \tilde{h}(\rho)\; \forall \rho$.
\end{theorem}
\proof The proof follows from Theorem~\ref{stable} and from all the arguments given above in this section,
in particular from $0~\leq \delta_j(\rho)~\leq~1, \forall \rho, j$. \endproof

\subsection{Numerical results}
\label{sec-sim}

We present in this section some numerical results.
We distinguish symmetric arrival passenger case where
the passenger arrival rates are the same to all the platforms, from
the asymmetric one where different passenger arrival rates are assumed.

With the parameters taken in Table~\ref{tab-param} above and in Table~\ref{tab-param2} below, we have $f_{\max} = 50$ train/h.,
$v = 41.18$ km/h, $w' = 26.61$ km/h, and $L = 17.294$ km. 
Let us notice here that $v_{\text{run}}$ given in Table~\ref{tab-param} corresponds to the free train speed during the inter-segment running times 
without including dwell times, while $v=L/\sum_j \underline{t}_j$ introduced in section~\ref{sec-ftd} corresponds to the train speed during
the whole travel time in the line, including train dwell times.

\begin{table}[htbp]
\centering
\caption{Additional parameters of the line considered.}
\begin{tabular}{|l||l|}
  \hline
  Passenger arrival rate (symmetric case) $\lambda$ & 1 pass./sec. (one platform)\\
  \hline  
  Passenger upload rate $\alpha$ & 30 passengers/s \\
  \hline
\end{tabular}
\label{tab-param2}
\end{table}

\begin{figure}[thbp]
\centering
  \includegraphics[scale=0.18]{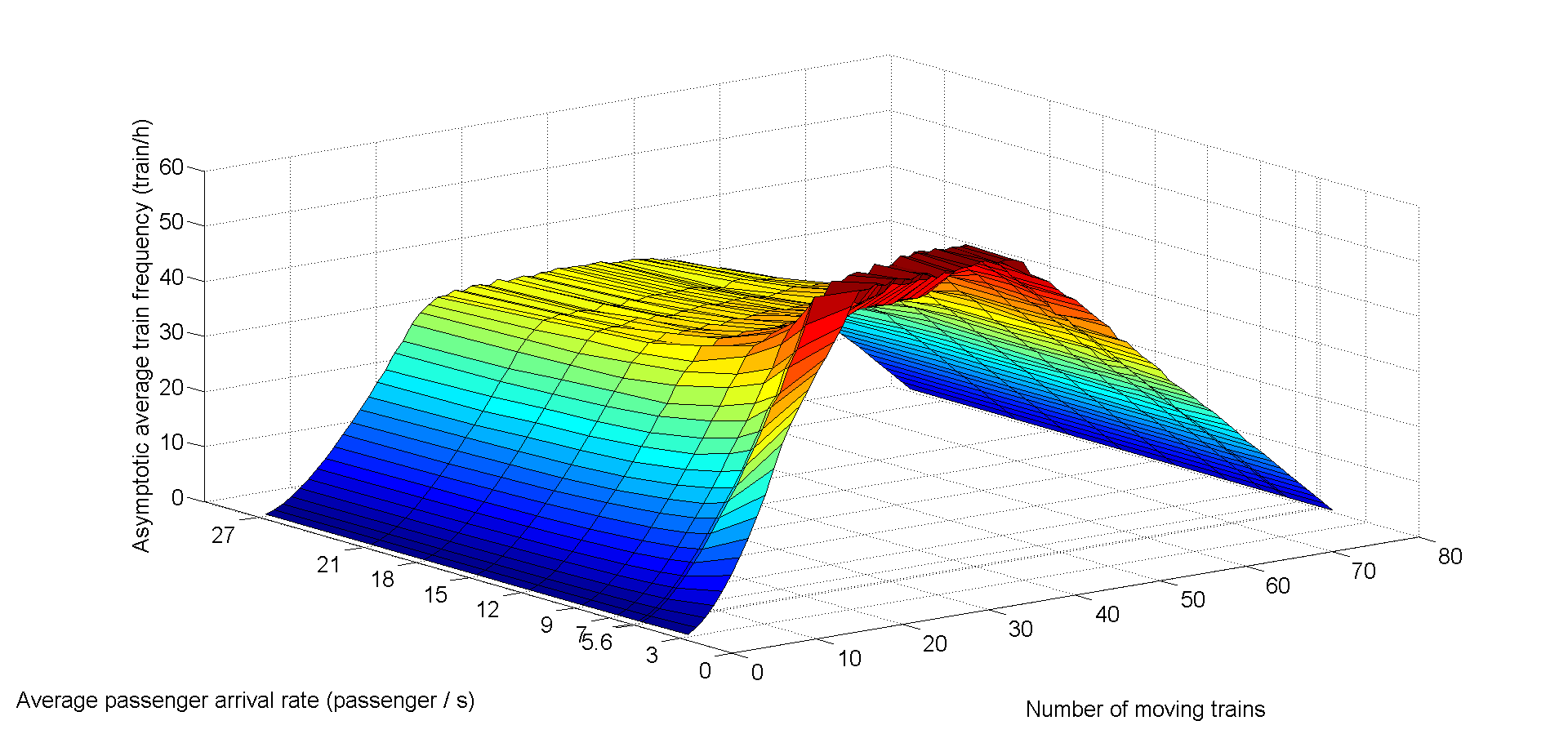} 
  \caption{Asymptotic average train frequency as a function of the number of running trains and of the average passenger arrival rates to the platforms
         (symmetric passenger arrival).}
\label{tab_2}
\end{figure}

\subsubsection*{Symmetric arrival passenger rate case}

In the symmetric arrival passenger rate case, 
the passenger arrival rate is assumed to be the same for all the platforms of the metro line.
The rate is varied in order to derive its effect on the train dynamics and on the physics of traffic.
The results are given in the left side of Table~\ref{tab_1} and in Figure~\ref{tab_2}.
The left side of Table~\ref{tab_1} shows the increasing of the train time-headways (degradation of the train frequencies) due to increases 
in the passenger arrival rates.
Figure~\ref{tab_2} gives a three dimensional illustration of this effect.

As shown in Theorem~\ref{th_stab}, the control law proposed here guarantees train dynamic stability, for every level of passenger demand.
For a passenger arrival rate less than or equal to the maximum supply of the line (maximum
passenger flow that can be served by the line), the control law guarantees the existence
of a number of trains, from which the train frequency is not affected by the passenger demand.

\subsubsection*{Asymmetric arrival passenger rate case}

\begin{figure}[htpb]
      \centering
      \includegraphics[scale=0.22]{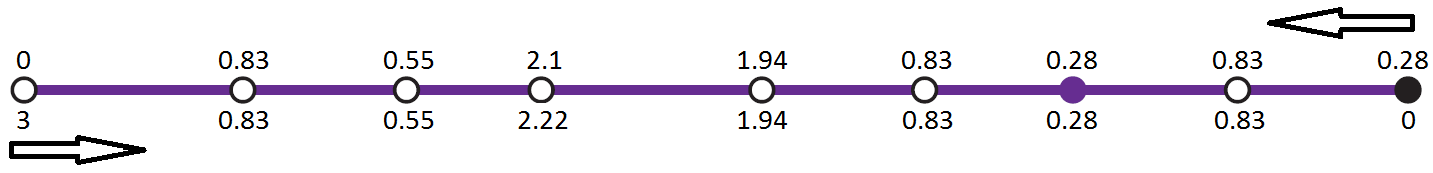}
      \caption{Average arrival rates $\lambda_j$ for every platform $j$, in passenger by second. The mean of those rates is 1.}
      \label{arrival-rates}
\end{figure}

Let us now consider the asymmetric passenger arrival rates on platforms.
We consider the arbitrary distribution of the average arrival rates of passengers shown in Figure~\ref{arrival-rates}.
The mean of those arrival rates over all the platforms is intentionally fixed to 1 passenger by second.

In the right-side figures of Table~\ref{tab_1}, we show the dependence of the asymptotic average train time-headway and frequency
on the intensity of passenger arrival rates.
We see that, although the mean of the average passenger arrival rates for both the symmetric and asymmetric distributions of those rates, 
is the same (equal to 1), the asymptotic train time-headways and frequencies differ.
We see on the right side figures of Table~\ref{tab_1} that the maximum of the average passenger arrival rates, 
over all the platforms counts more than the mean of those rates.
In the symmetric arrival distribution case, the maximum is equal to the mean and is~1 passenger by second.
In the asymmetric arrival distribution case, the maximum average arrival rate is~3; see Figure~\ref{arrival-rates}.
For example, the curves associated to $c=3$ in the right side of Table~\ref{tab_1} match better the curves
associated to the average arrival rate $c \max_j \lambda_j = 3\times 3 = 9$ in the left side of Table~\ref{tab_1}.
 
Finally, the figures of Table~\ref{tab_1} give the effect of the travel demand on the train dynamics, and on the train capacity of metro lines.
They permit also the determination of the minimum train time-headway (maximum frequency) of a given metro line, with known or predictable travel demand,
as well as the optimal number of running trains on the line, according to the travel demand.
Moreover, the changing of those indicators with respect to the number of running trains, is made available.
All this information can be used for planning, optimization, and real-time control of railway traffic in metro lines.

\section{CONCLUSIONS}

We proposed a dynamic programming-based approach for modeling and control of the train dynamics in metro lines.
The Max-plus linear model permits an analytical derivation of the traffic phases of the train dynamics in the case
where the passenger demand is not taken into account.
Although this first model is not realistic, it allows a good comprehension of the train dynamics.
Basing on the conclusions of the Max-plus linear model, we proposed an extension to a stochastic
dynamic programming model, where the passenger arrivals are taken into account.
The models permit the understanding of the effect of the passenger arrival demand on the train dwell times at platforms, and by that on the 
whole dynamics of the trains.
The perspectives in this direction of research are many.
First, the derivation of analytic formulas for the asymptotic average train frequency, dwell time and safe separation time 
for the stochastic dynamic programming model,
would bring a better comprehension of the traffic control model.
Second, the passenger demand being modeled here through average passenger arrival rates,
a dynamic model of the stock of passengers on the platforms an in the trains would improve the traffic dynamics, particularly by
taking into account the train and platform capacity limits.
Another direction of research is to extend the approach to metro lines with junctions.
Finally, other control parameters such as the train running times (speed profiles) can be considered in addition to train dwell times at platforms.

\appendices

\section{Proof of Theorem~\ref{th-dps2}}

\label{prf-dps2}

We show here that the dynamic system $x^k = \mathbf{f}(x^k,x^{k-1},\ldots,x^{k-m+1})$
is equivalent to another dynamic system $z^k = \mathbf{h}(z^{k-1})$, where $\mathbf{h}$ is built from $\mathbf{f}$.
Two steps are needed for the proof.

\textbf{Eliminate the dependence of $x^k$ on $x^{k-2}$, $x^{k-3}$, $\ldots$, $x^{k-m+1}$.}
    This is done by the state augmentation technique.
    In order to eliminate the dependence of $x_i^k$ on $x_j^{k-l}$, for every $1\leq i,j \leq n$, and for every $2\leq l\leq m-1$
    we add $l-1$ new variables $z_1,z_2, \ldots, z_{l-1}$ and the following $l$ dynamic equations to the dynamic system.
    $x_i^k = z_1^{k-1}, z_1^k = z_2^{k-1}, \ldots, z_{l-2}^k = z_{l-1}^{k-1}, z_{l-1}^k = x_j^{k-1}$.    
    By doing that, the new state vector of the dynamic system will only have terms of the zero and first orders.
    We notice here that the new dynamic equations added to the system satisfy the properties of additive homogeneity
    and monotonicity.
    Moreover, the connectivity of the new graph associated to the new dynamic system is preserved.
    That is, the graph associated to new system is strongly connected if and only if the graph associated to the original system is so.
    Indeed, the state augmentation done here only replaces arcs $(i \to j)$ in the original graph by paths 
    $(i \to i_1 \to i_2 \to \ldots \to i_{l-1} \to j)$ in the new graph, where $i_1, i_2, \ldots, i_{l-1}$ are
    intermediate nodes associated to the new variables $z_1, z_2, \ldots, z_{l-1}$ introduced here.        
    
\textbf{Eliminate implicit terms (the dependence of $x^k$ on $x^k$).}
    This is done by defining an order of updating the $n$ components of $x^k$, in
    such a way that no implicit term appears. This is possible because $\mathcal G(\mathbf{f}_0)$ 
    is acyclic. After the state augmentation procedure, we obtain a dynamic system of the form $z^k = \mathbf{g}(z^k,z^{k-1})$,
    where $z\in\mathbb R^{n'}$, with $n'\geq n$.
    By replacing implicit terms by explicit ones, we obtain at the end of this procedure another map for the dynamic system.
    We denote this map by $\mathbf{h}$, and the dynamic system is then written $z^k = \mathbf{h}(z^{k-1})$.
    We need now to show that the properties of monotonicity, additive 1-homogeneity, and strong connectedness are
    preserved by passing from $\mathbf{g}$ to $\mathbf{h}$.
        
    \textit{Homogeneity}. Replacing implicit terms by their explicit expressions, is nothing but a composition
         of two one-dimensional maps. Then, since additive homogeneity is preserved by composition, and since
         $\mathbf{g}$ is additive homogeneous, we conclude that $\mathbf{h}$ is so.
         
    \textit{Monotonicity}. Similarly, monotonicity being preserved by composition, and $\mathbf{g}$ being monotone, we conclude
         that $\mathbf{h}$ is monotone.
         
    \textit{Strong connectedness}. Let us use the notation
         $\mathbf{g}_0$ for the map $\mathbf{g}_0(x) = \mathbf{g}(x,-\infty)$, the notation $\mathbf{g}_1$ for the map $\mathbf{g}_1(x)=\mathbf{g}(-\infty,x)$, and the notation 
         $\tilde{\mathbf{g}}$ for the map $\tilde{\mathbf{g}}(x) = \mathbf{g}(x,x)$.
         We have $\mathbf{g}_0 = \mathbf{f}_0$, and then $\mathcal G(\mathbf{g}_0)$ is acyclic. Then two cases are distinguished.

           - If the graph $\mathcal G(\mathbf{g}_1)$ is strongly connected, then the graph $\mathcal G(\mathbf{h})$ will also be
             strongly connected, since all the arcs of $\mathcal G(\mathbf{g}_1)$ remain in $\mathcal G(\mathbf{h})$.
             
           - If $\mathcal G(\mathbf{g}_1)$ is not strongly connected, then it contains at least two strongly connected
             components, and eventual \textit{sink} nodes. The latter are nodes with no leaving (or outgoing) arcs in $\mathcal G(\mathbf{g}_1)$.
             We will show below that all the strongly connected components of the graph $\mathcal G(\mathbf{g}_1)$ will be
             connected into a unique strongly connected component in $\mathcal G(\mathbf{h})$, except the sink nodes which
             remain sinks in $\mathcal G(\mathbf{h})$.
             
             We assume that the graph $\mathcal G(\mathbf{g}_1)$ has only two strongly connected components.
             The proof is then easily extensible by induction to the case with more than two strongly connected components.
             
             Since $\mathcal G(\tilde{\mathbf{g}})$ is strongly connected, then there exists in $\mathcal G(\tilde{\mathbf{g}})$ at least one arc 
             going from a node in component 1 to a node in component 2, and at least one arc going from a node in
             component 2 to a node in component 1. One of these two arcs, or both, belong necessarily to $\mathcal G(\mathbf{g}_0)$,
             since $\mathcal G(\mathbf{g}_1)$ is not strongly connected.              
             We assume, without loss of generality, that such arc (the one or one of the two belonging to $\mathcal G(\mathbf{g}_0)$)
             is going from component 1 to component 2, and denote $n_1$ and $n_2$ its tail and head nodes respectively.
             We will prove below the following assertion.
             
             (A) : there exists an arc in $\mathcal G(\mathbf{h})$ from component 1 to component 2.
             
             If component 1 includes more than one node, then, $n_1$ has at least one incoming arc in $\mathcal G(\mathbf{g}_1)$,
             since component 1 is strongly connected in that graph.
             Let us denote $n_3$ the tail node of this incoming arc. 
             The arc $n_1 \to n_2$ interprets the fact that $z_{n_2}^k$ depends on $z_{n_1}^k$ in an implicit way.
             The arc $n_3 \to n_1$ interprets the fact that $z_{n_1}^k$ depends on $z_{n_3}^{k-1}$ in an explicit way.                          
             Therefore, $z_{n_2}^k$ depends on $z_{n_3}^{k-1}$ in an explicit way (by replacing $z_{n_1}^k$).
             That is, by elimination of the implicit term associated to the arc $n_1 \to n_2$, we obtain an (explicit) arc in $\mathcal G(\mathbf{h})$
             from node $n_3$ to node $n_2$, i.e. from component 1 to component 2. Hence (A) holds.
                          
             If component 1 includes a unique node $n_1$, then,              
             since $n_1$ is not a sink node, it has necessarily at least one outgoing arc in $\mathcal G(\mathbf{g}_1)$.
             If the outgoing arc from node $n_1$ goes to component~2, then we have shown that there exists an arc in 
             $\mathcal G(\mathbf{g}_1)$ from component 1 to component~2. This arc remain valid in $\mathcal G(\mathbf{h})$. Hence (A) holds.
             Otherwise, the outgoing arc from node $n_1$ is a loop arc.
             therefore, the arc is also an incoming arc to node $n_1$. 
             Then by denoting by $n_3$ the tail node of this incoming arc, and follow the same reasoning as above, 
             we obtain an (explicit) arc in $\mathcal G(\mathbf{h})$
             from node $n_3$ to node $n_2$, i.e. from component 1 to component 2. Hence (A) holds.
                                                                                         
             Now, if the arc going from component 2 to component 1 in $\mathcal G(\tilde{\mathbf{g}})$ is in $\mathcal G(\mathbf{g}_0)$, then
             we can show by the same way that it will be replaced by another arc from component 2 to component 1 in $\mathcal G(\mathbf{h})$.
             If the arc going from component 2 to component 1 in $\mathcal G(\tilde{\mathbf{g}})$ is in $\mathcal G(\mathbf{g}_1)$, then we will
             have the same arc in $\mathcal G(\mathbf{h})$. Therefore, there will be necessarily an arc from component 2
             to component 1 in $\mathcal G(\mathbf{h})$.
             
             By consequent, components 1 and 2 are connected into one strongly connected component in $\mathcal G(\mathbf{h})$.             
             
             Now, back to the sink nodes. Those nodes will remain sinks in the graph $\mathcal G(\mathbf{h})$ after replacing all
             implicit terms. They do not have any effect on the dynamics, since the associated variables only undergo
             the dynamics of the other variables.             

\addtolength{\textheight}{-17.2cm}   


\end{document}